\documentclass[11pt, a4paper, twoside]{article}
\usepackage{amsfonts}
\usepackage{mathrsfs}
\usepackage{amsmath}
\usepackage{amssymb}
\usepackage{fancyhdr}
\usepackage{hyperref}
\usepackage{color}
\numberwithin{equation}{section}

\allowdisplaybreaks

\begin{document}

\fancyhf{}

\fancyhead[EC]{Q. Xue}

\fancyhead[EL]{\thepage}

\fancyhead[OC]{Weighted estimates for the iterated Commutators of Multilinear Operators}

\fancyhead[OR]{\thepage}

\renewcommand{\headrulewidth}{0pt}
\renewcommand{\thefootnote}{\fnsymbol {footnote}}

\title{\textbf{Weighted Estimates for the iterated Commutators of Multilinear Maximal and Fractional Type Operators}
\footnotetext {{}{2000 \emph{Mathematics Subject
 Classification}: 42B20, 42B25}}
\footnotetext {{}\emph{Key words and phrases}: Multilinar Calder\'{o}n-Zygmund
operators, Maximal operators, Multilinear
fractional type operators, Commutators, multiple weights $A_{\vec{p}}$ and $A_{(\vec {p},q)}.$}} \setcounter{footnote}{0}
\author{
Qingying Xue \footnote {The author was supported
partly by NSFC (Grant No.10701010), NSFC (Key program Grant
No.10931001), PCSIRT of China, Beijing Natural
Science Foundation (Grant: 1102023).}\\
\hspace {0.1cm}}
\date{}
\maketitle

\begin{abstract}In this paper, the following iterated commutators $T_{*,\Pi b}$ of maximal operator for multilinear singular integral operators and $I_{\alpha, \Pi b}$ of multilinear fractional integral operator are introduced and studied
$$\aligned
T_{*,\Pi b}(\vec{f})(x)&=\sup_{\delta>0}\bigg|[b_1,[b_2,\cdots[b_{m-1},[b_m,T_\delta]_m]_{m-1}\cdots]_2]_1 (\vec{f})(x)\bigg|,
\endaligned
$$

$$\aligned
I_{\alpha, \Pi b}(\vec{f})(x)&=[b_1,[b_2,\cdots[b_{m-1},[b_m,I_\alpha]_m]_{m-1}\cdots]_2]_1 (\vec{f})(x),
\endaligned
$$
where $T_\delta$ are the smooth truncations of the multilinear singular integral operators and $I_{\alpha}$ is the multilinear fractional integral operator, $b_i\in BMO$ for $i=1,...,m$ and $\vec {f}=(f_1,...,f_m)$.

Weighted strong and $L(\log
L)$ type end-point estimates for the above iterated commutators associated with two class of multiple weights  $A_{\vec{p}}$ and $A_{(\vec{p}, q)}$ are obtained, respectively.\end{abstract}

\newtheorem{theorem}{Theorem}[section]
\newtheorem{definition}{Definition}[section]
\newtheorem{lemma}{Lemma}[section]
\newtheorem{proposition}{Proposition}[section]
\newtheorem{corollary}{Corollary}[section]
\newtheorem{remark}{Remark}[section]

\section[Introduction]{Introduction}

The multilinear Calder\'{o}n-Zygmund theory is a natural
generalization of linear case. Many authors were interested in these topics
(\cite{on commutators},
 \cite{operators}, \cite{Polynomial}, \cite{kenig frac},
\cite{multi C-Z}, \cite {DTT}, \cite{new maximal
function}, \cite{moen}, \cite{CX}, \cite{LXY}, \cite{PPTT}, \cite{GLPT} and \cite{BGP}). So we first recall the
definition and some results of multilinear
Calder\'{o}n-Zygmund operators as well as the corresponding multilinear maximal operators and fractional type operators.\\
\begin{definition} [Multilinear Calder\'{o}n-Zygmund
operators]  Let T be a  Multilinear operator initially defined on
the m-fold product of Schwartz spaces and taking values in the space
of tempered distributions,
\begin{equation*}
T:\mathscr{S}(\mathbb{R}^{n})\times\cdots\times\mathscr{S}(\mathbb{R}^{n})
\longrightarrow\mathscr{S}'(\mathbb{R}^{n}).
\end{equation*}
Following [6], we say that $T$ is an $m$-linear
Calder\'{o}n-Zygmund operator if for some $1\leq q_{j}<\infty$, it
extends to a bounded multilinear operator from
$L^{q_{1}}\times\cdots\times L^{q_{m}}$ to $L^{q}$, where
$\frac{1}{q}=\frac{1}{q_{1}}+\cdots+\frac{1}{q_{m}}$, and if there
exists a function $K$, defined off the diagonal $x=y_{1}=\cdots=y_{m}$
in $(\mathbb{R}^{n})^{m+1}$, satisfying
\begin{equation*}
T(f_{1},\cdots,f_{m})(x)
=\int_{(\mathbb{R}^{n})^{m}}K(x,y_{1},\cdots,y_{m})f_{1}(y_{1})\cdots
f_{m}(y_{m})dy_{1}\cdots dy_{m}
\end{equation*}
for all $x\not\in\bigcap_{j=1}^{m}\operatorname{supp}f_{j};$
\begin{equation}
|K(y_{0},y_{1},\cdots,y_{m})|\leq\frac{A}
{(\sum_{k,l=0}^{m}|y_{k}-y_{l}|)^{mn}};
\end{equation}
and
\begin{equation}
|K(y_{0},\cdots,y_{j},\cdots,y_{m})-K(y_{0},\cdots,y_{j}',\cdots,y_{m})|
\leq\frac{A|y_{j}-y_{j}'|^{\varepsilon}}
{(\sum_{k,l=0}^{m}|y_{k}-y_{l}|)^{mn+\varepsilon}},
\end{equation}
for some $\varepsilon>0 $ and all 0$\leq j\leq m$, whenever
$|y_{j}-y_{j}'|\leq\frac{1}{2}\max_{0\leq k \leq m}|y_{j}-y_{k}|.$\\

The maximal multilinear singular integral operator was defined by
$$T_{*}(\vec{f})(x)=\sup_{\delta>0}|T_{\delta}(f_1,\cdots,f_m)(x)|,$$
where $T_{\delta}$ are the smooth truncations of T given by
$$T_{\delta}(f_1,\cdots,f_m)(x)=\int_{|x-y_1|^2+\dots+|x-y_m|^2>{\delta}^2}
K(x,y_1,\cdots,y_m)f_1(y_1)\cdots f_m(y_m)d\vec{y}.$$ \\
Here, $d\vec{y}=dy_1\cdots dy_m$.\end{definition}
As is pointed in \cite{GT1}, $T_*(\vec{f})(x)$ is pointwise well-defined when $f_j\in L^{q_j}(\mathbb{R}^n)$ with $1\le q_j\le \infty.$

The study of the multilinear singular integral operator and its  maximal operator has a long history. For maximal multilinear operator $T_{*}$, one can see for example,
 \cite{GT1}, \cite{on Hardy spaces}, \cite {LXY} and \cite {C} for more details. We list some results for
${T_{*}}$ as follows:

\textbf{Theorem A}{(\cite{GT1})} Let $1\leq q_i <\infty$, and $q$ be such
that $\frac{1}{q}=\frac{1}{q_1}+\cdots+\frac{1}{q_m}$, and $\omega
\in A_{q_1}\cap \cdots \cap A_{q_m}$. Let T be an m-linear
Calder\'{o}n-Zygmund operator. Then there exists a constant
$C_{q,n}<\infty$ so that for all $\vec{f}=(f_1,\cdots,f_m)$
satisfying $$\|{T_*}(\vec{f})\|_{L^{q}_{\omega}}\leq
C_{n,q}(A+W)\prod_{i=1}^m\|f_i\|_{L^{q_i}_{\omega}},$$
where W is the norm of T in the mapping T: $L^1 \times \cdots \times L^1
\rightarrow L^{1/m,\infty}$.

\textbf{Theorem B }(\cite{C}) Assume that $\frac{1}{p_{1}}+\cdots+\frac{1}{p_{m}}=\frac{1}{p}$ and $\vec {w}\in A_{\vec{p}}$, then\begin{enumerate}
\item [(i)] If $1<p_1,...,p_m<\infty$, then $T_*$ is bounded from $L^{p_1}(w_1)\times \cdots \times L^{p_m}(w_m)$ to $L^p(\vec{\omega});                                $
\item [(ii)] If $1\le p_1,...,p_m<\infty$, then $T_*$ is bounded
from $L^{p_1}(w_1)\times \cdots \times L^{p_m}(w_m)$ to
$L^{p,\infty}(\vec{\omega}).$\end{enumerate}

Here, $A_{\vec{p}}$ is the multiple weights in the Definition 2.1 below. The boundedness of $T_*$ on Hardy spaces and weighted Hardy spaces were obtained in \cite{on Hardy spaces} and \cite{LXY1}.

Now, let's recall some definitions and background for the multilinear fractional type operators.

In 1992, Grafakos \cite{multi frac grafakos} first defined and
studied the multilinear maximal function and multilinear fractional
integral as follows
\begin{equation*}M_{\alpha}(\vec{f})(x) = \sup_{r > 0}\frac{1}{r^{n
- \alpha}}\int_{|y| < r}\Big|\prod_{i = 1}^mf_i(x -
\theta_iy)\Big|\,dy\end{equation*} and
\begin{equation*}I_{\alpha}(\vec{f})(x) =
\int_{\mathbb{R}^n}\frac{1}{{|y|}^{n - \alpha}}\prod_{i = 1}^mf_i(x
- \theta_iy)\,dy,\end{equation*}where $\theta_i$ $(i = 1, \cdots,
m)$ are fixed distinct and nonzero real numbers and $0 < \alpha <n$.
We note that, if we simply take $m=1$ and $\theta_i=1$, then
$M_\alpha$ and $I_\alpha$ are just the operators studied by
Muckenhoupt and Wheeden in \cite{MW}. In 1999, Kenig and Stein
\cite{kenig frac} considered another more general type of
multilinear fractional integral which was defined by
\begin{equation*}I_{\alpha, A}(\vec{f})(x) = \int_{{(\mathbb{R}^n)}^m}\frac{1}{{|(y_1, \cdots, y_m)|}^{mn - \alpha}}
\prod_{i = 1}^mf_i\big(\ell_i(y_1, \cdots, y_m,
x)\big)\,dy_i,\end{equation*} where $\ell_i$ is a linear combination
of $y_j$s and $x$ depending on the matrix $A$. They showed that
$I_{\alpha, A}$ was of strong type $(L^{p_1} \times \cdots \times
L^{p_m}, \ L^q)$
 and weak type $(L^{p_1} \times \cdots \times L^{p_m}, \ L^{q,\infty})$.
When $\ell_i(y_1,\cdots,y_m,x)=x-y_i$, we
denote this multilinear fractional type operator by
$I_{\alpha}$.

\indent For a long time, there is an open question (\cite{report}) in the multilinear operators theory. That is, the existence of multiple weights theory for
multilinear Calder\'{o}n-Zygmund operators and multilinear fractional integral operators. This was established in
\cite{new maximal function}, \cite{moen}, \cite{CX} and the multiple weights
$A_{\vec{p}}$ and $A_{(\vec{p}, q)}$ were constructed (see the definitions in section 2 below).

In \cite{new maximal function} and \cite{CX}, the following commutators of $T$ and $I _\alpha$ in the j-th entry were defined and studied, including weighted strong and weighted end-point $L(\log L)$ type estimates associated with $A_{\vec {p}}$ and $A_{(\vec{p}, q)}$ weights, respectively.

\begin{definition} [Commutators in the j-th entry](\cite{new maximal function}, \cite{CX}) Given a
collection of locally integrable functions
$\vec{b}=(b_1,\cdots,b_m)$, we define the commutators of
the m-linear Calder\'{o}n-Zygmund operator $T$ and fractional integral $I_\alpha $ to be
$$[\vec {b},T](\vec{f})=T_{\vec{b}}(f_1,\dots,f_m)=\sum_{j=1}^mT_{\vec{b}}^j(\vec{f}),\quad I_{\vec{b}, \alpha}(\vec{f})(x) = \sum_{i = 1}^mI_{\vec{b}, \alpha}^i(\vec{f})(x),$$
where each term is the commutator of $b_j$ and T in the j-th entry
of T, that is,
$$T_{\vec{b}}^j(\vec{f})=b_jT(f_1,\dots,f_j,\dots,f_m)-T(f_1,\dots,b_jf_j,\dots,f_m).$$
\begin{equation*}\end{equation*}
Also\begin{equation*}I_{\vec{b},
\alpha}^i(\vec{f})(x) = b_i(x)I_\alpha(f_1, \cdots, f_i, \cdots,
f_m)(x) - I_\alpha(f_1, \cdots, b_if_i, \cdots,
f_m)(x).\end{equation*}
\end{definition}
\indent
Recently, in \cite{PPTT}, the following iterated commutators of multilinear Calder¡äon-Zygmund operators and
pointwise multiplication with functions in BMO are defined and studied in products of Lebesgue
spaces, including strong type and weak end-point estimates with multiple $A_{\vec{p}}$ weights.
\begin{equation}\aligned
T_{\Pi
b}(\vec{f})(x)&=[b_1,[b_2,\cdots[b_{m-1},[b_m,T]_m]_{m-1}\cdots]_2]_1
(\vec{f})(x)\\&=\int_{{(\mathbb{R}^n)}^m}{\prod
_{j=1}^m(b_j(x)-b_j(y_j))K(x,y_1,...,y_m)} \prod_{i =
1}^{m}f_i(y_i)\,d\vec{y}.
\endaligned
\end{equation}
\indent Therefore, an open interesting question arises, can we establish the weighted strong and end-point estimates of the iterated commutators for the multilinear operator $T_*$ and $I_\alpha$? We note that, there is no results for the commutators of multilinear operator $T_*$ $(m\ge 2)$, even for the commutators of $T_*$ in the j-th entry.

\indent In this article, we give a positive answer to the above question, we study iterated commutators of maximal multilinear singular integral operator and multilinear fractional integral operators
defined by
\begin{equation}\aligned
T_{*,\Pi b}(\vec{f})(x)&=\sup_{\delta>0}\bigg|[b_1,[b_2,\cdots[b_{m-1},[b_m,T_\delta]_m]_{m-1}\cdots]_2]_1 (\vec{f})(x)\bigg|\\&=\sup_{\delta>0}\bigg|\int_{|x-y_1|^2+\dots+|x-y_m|^2>{\delta}^2}{\prod _{j=1}^m(b_j(x)-b_j(y_j))K(x,y_1,...,y_m)}
\prod_{i = 1}^{m}f_i(y_i)\,d\vec{y}\bigg|
\endaligned
\end{equation}
and
\begin{equation}\aligned
I_{\alpha, \Pi b}(\vec{f})(x)&=[b_1,[b_2,\cdots[b_{m-1},[b_m,I_\alpha]_m]_{m-1}\cdots]_2]_1 (\vec{f})(x)\\&=\int_{{(\mathbb{R}^n)}^m}\frac{1}{{|(x-y_1, \cdots, x-y_m)|}^{mn - \alpha}}\prod _{j=1}^m(b_j(x)-b_j(y_j))
\prod_{i = 1}^{m}f_i(y_i)\,d\vec{y}.
\endaligned
\end{equation}
\begin{remark} Note that, when $m=1$ in (1.3), this definition coincides with the linear commutator
$[b,T]f=bT(f)-T(bf)$ and $[b,I_\alpha]f=bI_\alpha(f)-I_\alpha(bf)$.
One classical result given by Coifman, Rochberg and Weiss \cite{CRW}
is that $[b,T]$ is $L^p$ bounded for $1<p<\infty$ when $b\in BMO$.
But $[b,T]$ fails to be an operator of weak type $(1,1)$, a
counterexample was given by C. P\'{e}rez and an alternative $L(\log
L)$ type result was obtained in \cite {pe}. In 1982, Chanillo proved
that the commutator of the fractional integral operator $[b,
I_\alpha]$ is bounded from $L^p$ into $L^q$ ($p> 1, 1/q = 1/p -
a/n$) when $b\in BMO$. In 2002, Ding, Lu and Zhang \cite{DLZ}
studied the continuity properties of fraction type operators. They
showed that $[b,I_\alpha]$ fails to be an operator of weak type
$(L^1,L^{n/(n-\alpha),\infty})$, counterexamples were given in
\cite{DLZ}, alternative $L(\log L)$ type estimates was obtained.
\end{remark}

We state our results as follows.

\begin{theorem}[Weighted strong bounds for $T_{*,\Pi b}$]\label{strong commutator}
Let $\vec{\omega}\in A_{\vec{p}}$, $\frac{1}{p}
= \frac{1}{p_1} + \cdots + \frac{1}{p_m}$ with $1<p_j<\infty$, $j=1,..,m$; and $\vec{b}\in (BMO)^m$. Then
 there is a constant $C > 0$ independent of $\vec{b}$ and $\vec{f}$ such
that\begin{equation}\label{weighted strong commutator}
 {\big\|T_{*,\Pi b}(\vec{f})\big\|}_{L^{p}({\nu_{\vec{\omega}}})}
\leqslant C {\prod_{j=1}^{m}{\big\|b_j\big\|}_{{BMO}}}\prod_{i =
1}^m{\big\|f_i\big\|}_{L^{p_i}(\omega_i)},\end{equation}where
$\vec {b}= (b_1,...,b_m)$.\end{theorem}

\begin{theorem}[Weighted end-point estimate for $T_{*,\Pi b}$]\label{end-point commutator}
Let $\vec{\omega} \in A_{(1, \cdots, 1))}$ and $\vec{b} \in {(BMO)}^m$.
 Then there exists a constant $C$ depending on $\vec{b}$ such that
\begin{equation}\label{commutator weak}{\nu_{\vec{\omega}}}
\bigg(\Big\{x\in \mathbb{R}^n :T_{*,\Pi b}(\vec{f})(x) > t^{{m}}\Big\}\bigg)
 \leqslant C
 {\bigg(\prod_{i = 1}^m\int_{\mathbb{R}^n}\Phi ^{(m)}\Big(\frac{|f_i(y_i)|}{t}\Big)\omega_i(y_i)\,dy_i\bigg)}^{\frac{1}{m}},\end{equation}
 \end{theorem}
 where $\Phi (t)=t(1+\log ^+t)$ and $\Phi ^{(m)}=\overbrace{\Phi \circ\cdots\circ\Phi}^{m} $.
\begin{remark}
If $m=1$, then weighted strong $L^p$ and weighted end-point $L(\log L)$ estimates for commutators of the classical linear operator $T_{*}$ were studied in \cite {zhang}.
\end{remark}
As for $I_{\alpha, \Pi b}$, we get
\begin{theorem}[Weighted strong bounds for $I_{\alpha, \Pi b}$]\label{strong commutator}
Let $0 < \alpha < mn$, $1 < p_1, \cdots, p_m < \infty$, $\frac{1}{p}
= \frac{1}{p_1} + \cdots + \frac{1}{p_m}$ and $\frac{1}{q} =
\frac{1}{p} - \frac{\alpha}{n}$. For $r > 1$
 with $0 < r \alpha < mn$, if \ $\vec{\omega}^{r}\in A_{(\frac{\vec{p}}{r},
 \frac{q}{r})}$, ${\nu_{\vec{\omega}}}^q \in A_\infty$ and $\vec{b} \in {(BMO)}^m$,
 there is a constant $C > 0$ independent of $\vec{b}$ such
that\begin{equation}\label{weighted strong commutator}
 {\big\|I_{\alpha, \Pi b}(\vec{f})\big\|}_{L^{q}({\nu_{\vec{\omega}}}^q)}
\leqslant C {\prod_{j=1}^{m}{\big\|b_j\big\|}_{{BMO}}}\prod_{i =
1}^m{\big\|f_i\big\|}_{L^{p_i}(\omega_i^{p_i})}.\end{equation}\end{theorem}

\begin{theorem}[Weighted end-point estimate for $I_{\alpha, \Pi b}$]\label{end-point commutator}
Let $0 < \alpha < mn$, $\vec{\omega} \in A_{((1, \cdots, 1),
\frac{n}{mn - \alpha})}$ and $\vec{b} \in {(BMO)}^m$.
 Then there exists a constant $C$ depending on $\vec{b}$, such that
\begin{equation}\aligned\label{commutator weak}{\nu_{\vec{\omega}}}&^{\frac{n}{mn -\alpha}}
\bigg(\Big\{ x\in \mathbb{R}^n : I_{\alpha, \Pi b}(\vec{f})(x) > t^{\frac{mn - \alpha}{n}}\Big\}\bigg)\\&
 \leqslant C
 {\bigg\{ \bigg[1 + \frac{\alpha}{mn}\log^+\bigg(\prod_{i = 1}^m\int_{\mathbb{R}^n}\Phi^{(m)}\Big(\frac{|f_i(y_i)|}{t}\Big)\,dy_i\bigg)\bigg]^m \prod_{j = 1}^m\int_{\mathbb{R}^n}\Phi^{(m)}\Big(\frac{|f_j(y_j)|}{t}\Big)\omega_j(y_j)\,dy_j\bigg\}}^{\frac{n}{mn - \alpha}}.\endaligned
 \end{equation}
 Moreover, if each $0<\alpha_j<n$, we obtain
 \begin{equation}\aligned\label{commutator weak}{\nu_{\vec{\omega}}}&^{\frac{n}{mn -\alpha}}
\bigg(\Big\{ x\in \mathbb{R}^n : I_{\alpha, \Pi b}(\vec{f})(x) > t^{\frac{mn - \alpha}{n}}\Big\}\bigg)\\&
 \leqslant C
 {\bigg\{\prod_{j = 1}^m \bigg[1 + \frac{\alpha_j}{n}\log^+\bigg(\prod_{i = 1}^m\int_{\mathbb{R}^n}\Phi^{(m)}\Big(\frac{|f_i(y_i)|}{t}\Big)\,dy_i\bigg)\bigg]\int_{\mathbb{R}^n}\Phi^{(m)}\Big(\frac{|f_j(y_j)|}{t}\Big)
 \omega_j(y_j)\,dy_j\bigg\}}^{\frac{n}{mn - \alpha}},\endaligned
 \end{equation}
 where $\Phi (t)$ and $\Phi ^{(m)}$ are the same as in Theorem 1.2.
 \end{theorem}
As a corollary of Theorem 1.3 and Theorem 1.4, we can obtain similar results for the commutators of the multilinear fractional maximal operator. Let's first give its definition. Suppose each $f_i$ $(i = 1, \cdots,
m)$ is locally integrable on $\mathbb{R}^n$. Then for any $x\in
\mathbb{R}^n$, we define the multilinear fractional maximal operator and its commutators
by
\begin{equation*}\mathcal{M}_\alpha(\vec{f})(x) =
\sup_{Q}{|Q|}^{\frac{\alpha}{n}}\prod_{i =
1}^m\frac{1}{|Q|}\int_{Q}|f_i(y_i)|\,dy_i,\end{equation*}
and
\begin{equation*}\mathcal{M}_{{\alpha,\Pi b}}(\vec{f})(x) =
\sup_{Q}{|Q|}^{\frac{\alpha}{n}}\prod_{i =
1}^m\frac{1}{|Q|}\int_{Q}|b_i(x)-b_i(y_i)||f_i(y_i)|\,dy_i,\end{equation*}
where the
supremum is taken over all cubes $Q$ containing $x$ in
$\mathbb{R}^n$ with the sides parallel to the axes.
\begin{corollary}
 Let $\alpha$, $b_i$, $\vec{\omega}$, $p_i, q$ be the same as in Theorem 1.3-1.4, then Theorem 1.3-1.4 still hold for $\mathcal{M}_{{\alpha,\Pi b}}.$ \end{corollary}

The article is organized as follows. In section \ref{De}, we prepare
some definitions and lemmas. Some propositions will be listed and
proved in section \ref{def}, including the main Proposition 3.1.
Then, we give the proof of Theorem 1.1-1.3. Section \ref{commutators
2} will be devoted to the study of the end-point $L(\log L)$ type
estimates for the iterated commutators of multilinear fractional
type operators.

\section[Definitions and some lemmas]{Definitions and some lemmas}\label{De}

Let us recall the definitions of $A_{\vec{p}}$ and $A_{(\vec{p}, q)}$ weights.\\
For $m$-exponents $p_{1},\cdots,p_{m}$, we will often write $p$ for
the number given by
$\frac{1}{p}=\frac{1}{p_{1}}+\cdots+\frac{1}{p_{m}}$, and
$\vec{p}$ for the vector $\vec{p}=(p_{1},\cdots,p_{m})$.\\
\begin{definition}[Multiple $A_{\vec{p}}$ weights](\cite{new maximal function})
Let $1\leq p_{1},\cdots,p_{m}<\infty$. Given
$\vec{\omega}=(\omega_{1},\cdots,\omega_{m})$, set
\begin{equation*}
\nu_{\vec\omega}=\prod_{j=1}^{m}\omega_{j}^{p/p_{j}}.
\end{equation*}
We say that $\vec{\omega}$ satisfies the $A_{\vec{p}}$ condition if
\begin{equation}
\sup_Q{\bigg(\frac{1}{|Q|}\int_Q\prod_{i =
1}^m{\omega_i}^{\frac{p}{p_i}}\,\bigg)}^{\frac{1}{p}}
 \prod_{i = 1}^m{\bigg(\frac{1}{|Q|}\int_Q{\omega_i}^{1 -
p'_i}\,\bigg)}^{\frac{1}{p'_i}} < \infty.
\end{equation}
When $p_{j}=1,\bigg(\frac{1}{|Q|}\int_Q{\omega_i}^{1 - p'_i}\,\bigg)
^{\frac{1}{p'_i}}$ is understood as $(\inf_Q \omega_{j})^{-1}$.\\
\end{definition}

\begin{definition}[Multiple $A_{(\vec{p}, q)}$ weights](\cite{CX}, \cite{moen}) Let\ $1 \leqslant p_1, \cdots, p_m < \infty$, $\frac{1}{p} =
\frac{1}{p_1} + \cdots + \frac{1}{p_m}$, and $q > 0$. Suppose that
$\vec{\omega} = (\omega_1, \cdots, \omega_m)$ and each $\omega_i$
$(i = 1, \cdots, m)$ is a nonnegative function on $\mathbb{R}^n$. We
say that $\vec{\omega}\in A_{(\vec{p}, q)}$ if it
satisfies\begin{equation}\label{weight}\sup_Q{\bigg(\frac{1}{|Q|}\int_Q{\nu_{\vec{\omega}}}^q\,\bigg)}^{\frac{1}{q}}
\prod_{i =
1}^m{\bigg(\frac{1}{|Q|}\int_Q{\omega_i}^{-p'_i}\,\bigg)}^{\frac{1}{p'_i}}
< \infty,\end{equation}where $\nu_{\vec{\omega}} = \prod_{i =
1}^m\omega_i$. If $p_i = 1$,
$(\frac{1}{|Q|}\int_Q\omega_i^{-p'_i})^{\frac{1}{p'_i}}\,$ is
understood as $(\inf_Q\omega_i)^{-1}$.
\end{definition}

\begin{remark}
\indent In particular, when $m=1$, we note that $A_{\vec{p}}$ will
be degenerated to the classical $A_{p}$ weights. Moreover, if m=1 and
$p_i=1$, then  this class of weights coincide with the classical
$A_1$ weights. Also, when $m=1$, we note that $A_{(\vec{p}, q)}$ will be degenerated to the classical $A_{(p,
q)}$ weights, where the latter was defined in 1974 by B. Muckenhoupt and R. Wheeden \cite{MW}. We will refer to (1.4) and (1.5) as the multilinear
$A_{\vec{p}}$ condition and $A_{(\vec {p},
q)}$ condition.
\end{remark}
We need the following $L(\log)L$ type multilinear maximal fractional operators
\begin{definition}For any $\vec{f} = (f_1, \cdots, f_m)$ and $0 < \alpha <m n$ with
$\sum_{i = 1}^m\alpha_i = \alpha$, two multilinear fractional
$L(\log L)$ type maximal operators are defined as
\begin{equation*}\mathcal{M}_{L(\log L), \alpha}^j(\vec{f})(x)
 = \sup_{Q\ni x}{|Q|}^{\frac{\alpha}{n}}{\|f_j\|}_{L(\log L), Q}
\prod_{i \neq j}\frac{1}{{|Q|}}\int_Q|f_i|\,\end{equation*}and
\begin{equation*}\mathcal{M}_{L(\log L), \alpha}(\vec{f})(x)
 = \sup_{Q\ni x}{|Q|}^{\frac{\alpha}{n}}\prod_{i = 1}^m{\|f_i\|}_{L(\log L),
Q},
\end{equation*}respectively. If $\alpha=0$, for simply, we denote $\mathcal{M}_{L(\log L), 0}=\mathcal{M}_{L(\log L)}$ and $\mathcal{M}_{L(\log L), 0}^j=\mathcal{M}_{L(\log L)}^j$ \end{definition}

We prepare some lemmas which will be used later. The following H\"{o}lder's inequality on Orlicz
spaces can be seen in \cite[p. 58]{orlicz}.

\begin{lemma}[Generalized H\"{o}lder's inequality](\cite{orlicz})  Let $\phi(t) = t(1 + {\log}^+t)$ and
 $\psi(t) = e^t - 1$ and suppose
that\begin{eqnarray*}{\|f\|}_{\phi} & \triangleq &
 \inf\bigg\{\lambda > 0 : \int_{\mathbb{R}^n}\phi\Big(\frac{|f(x)|}{\lambda}\Big)\,d\mu \leqslant 1\bigg\}
 < \infty\\
{\|g\|}_{\psi} & \triangleq & \inf\bigg\{\lambda > 0 :
\int_{\mathbb{R}^n}\psi\Big(\frac{|g(x)|}{\lambda}\Big)\,d\mu
\leqslant 1\bigg\} < \infty\end{eqnarray*}with respect to some
measure $\mu$, then for any cube
$Q$\begin{equation}\label{generalized
holder}\frac{1}{|Q|}\int_Q|fg|\, \leqslant 2 {\|f\|}_{L(\log L),
Q}{\|g\|}_{\exp L, Q}.\end{equation}\end{lemma}

Some other inequalities are also necessary.

\begin{lemma} (\cite {CX}) Suppose that $r > 1$ and $b \in BMO$, then for any $f$
satisfing the condition of generalized H\"{o}lder's inequality there
is a $C
> 0$ independent of $\vec{f}$ and $b$ such that
\begin{eqnarray}\label{s1}\frac{1}{|Q|}\int_Q|f|\, & \leqslant & C {\|f\|}_{L(\log L),
Q};\\
\label{s2}{\|f\|}_{L(\log L), Q} & \leqslant & C
{\bigg(\frac{1}{|Q|}\int_Q{|f|}^r\,\bigg)}^{\frac{1}{r}};\\
\label{s3}\frac{1}{|Q|}\int_Q|(b - b_Q)f|\, & \leqslant & C
{\|b\|}_{BMO} {\|f\|}_{L(\log L),
Q};\\
\label{s4}{\bigg(\sup_Q\frac{1}{|Q|}\int_Q{|b - b_Q|}^{r -
1}\,\bigg)}^{\frac{1}{r - 1}} & \leqslant & C
{\|b\|}_{BMO}.\end{eqnarray}\end{lemma}

We need Kolmogorov's
inequalities in the following lemma, which are necessary tools for some estimates.
\begin{lemma}[Kolmogorov's inequality](\cite {new maximal function}, \cite[p. 485]{weighted norm inequ})
\begin{enumerate}
\item [(a)] Suppose $0<p<q<\infty$, then
\begin{equation}\|f\|_{L^p(Q, \frac {dx}{|Q|})}\le C \|f\|_{L^{q,\infty}(Q, \frac {dx}{|Q|})};\end{equation}
\item [(b)] Suppose that $0 < \alpha < n$ and $p, q > 0$ satisfying
$\frac{1}{q} = \frac{1}{p} - \frac{\alpha}{n}$. Then for any
measurable function $f$ and cube $Q$,
\begin{equation}\label{kolmogorov}{\bigg(\int_{Q}{|f|}^p\,\bigg)}^{\frac{1}{p}} \leqslant
{\bigg(\frac{q}{q - p}\bigg)}^{\frac{1}{p}} {\big|Q\big|}^{\frac{\alpha}{n}} {\big\|f\big\|}_{L^{q, \infty}(Q)}.\end{equation}.\end{enumerate}
\end{lemma}
To prove Theorem 1.4, we also need the following known results,
\begin{lemma}[Weighted estimates for $\mathcal{M}_{\alpha}$ and $I_\alpha$](\cite{moen}, \cite{CX})\label{weak maxi}
 Let $0 < \alpha < mn$, $1 \leqslant p_1, \cdots, p_m < \infty$,
$\frac{1}{p} = \frac{1}{p_1} + \cdots + \frac{1}{p_m}$ and
$\frac{1}{q} = \frac{1}{p} - \frac{\alpha}{n}$. Then for
$\vec{\omega}\in A_{(\vec{p}, q)}$ there is a constant $C > 0$
independent of $\vec{f}$ such that
\begin{equation}\label{Weighted Weak Boundedness of the Maximal Operator}
{\big\|\mathcal{M}_{\alpha}(\vec{f})\big\|}_{L^{q,
\infty}({\nu_{\vec{\omega}}}^{q})} \leqslant C \prod_{i =
1}^m{\big\|f_i\big\|}_{L^{p_i}(\omega_i^{p_i})};\end{equation}
\begin{equation}\label{Weighted Weak Boundedness of the fractional
integral}{\big\|I_{\alpha}(\vec{f})\big\|}_{L^{q,
\infty}({\nu_{\vec{\omega}}}^{q})} \leqslant C \prod_{i =
1}^m{\big\|f_i\big\|}_{L^{p_i}(\omega_i^{p_i})}.\end{equation}
\end{lemma}
\section[Proof of Thm 1.1-1.3]{Proof of Theorem 1.1-1.3}\label{def}

\vspace{0.5cm}

To begin with, we prepare one proposition which plays important role in the proof of our theorems. The basic idea is to control the iterated commutators of $T_*$ by another two operators.

Let $u, v\in C^\infty([0,\infty))$ such that $|u'(t)|\le Ct^{-1}$, $|v'(t)|\le Ct^{-1}$ and satisfy
$$\chi_{[2,\infty)}(t)\le u(t)\le \chi_{[1,\infty)}(t), \quad \chi_{[1,2]}(t)\le v(t)\le \chi_{[1/2,3]}(t).$$
We define the maximal operators
$$U^*(\vec{f})(x)=\sup _{\eta>0}\bigg|\int_{({\mathbb{R}}^n)^m}K(x,y_1,...,y_m)u(\sqrt{{|x-y_1|+,...,+|x-y_m|}}/{\eta})\prod_{i=1}^m
f_{i}(y_{i})d\vec{y}\bigg|,
$$
$$V^*(\vec{f})(x)=\sup _{\eta>0}\bigg|\int_{({\mathbb{R}}^n)^m}K(x,y_1,...,y_m)v(\sqrt{{|x-y_1|+,...,+|x-y_m|}}/{\eta})
\prod_{i=1}^mf_{i}(y_{i})d\vec{y}\bigg|.
$$
For simplicity, we denote $K_{u, \eta}(x,y_1,...,y_m)=K(x,y_1,...,y_m)u(\sqrt{{|x-y_1|+,...,+|x-y_m|}}/{\eta})$, $K_{v, \eta}(x,y_1,...,y_m)=K(x,y_1,...,y_m)v(\sqrt{{|x-y_1|+,...,+|x-y_m|}}/{\eta})$ and
$$U_\eta (\vec {f})=\int_{({\mathbb{R}}^n)^m}K_{u, \eta}(x,y_1,...,y_m)\prod_{i=1}^mf_{i}(y_{i})d\vec{y}$$ and
$$V_\eta (\vec {f})=\int_{({\mathbb{R}}^n)^m}K_{v, \eta}(x,y_1,...,y_m)\prod_{i=1}^mf_{i}(y_{i})d\vec{y}.$$
It is easy to see that $T_*(\vec{f})\le U^*(\vec{f})(x)+V^*(\vec{f})(x)$. Moreover, $T_{*,\Pi b}(\vec{f})\le U_{\Pi b}^*(\vec{f})(x)+V_{\Pi b}^*(\vec{f})(x)$,
where
$$\aligned
U_{\Pi b}^*(\vec{f})(x)&=\sup_{\eta>0}\bigg|[b_1,[b_2,\cdots[b_{m-1},[b_m,U_\eta ]_m]_{m-1}\cdots]_2]_1 (\vec{f})(x)\bigg|
\\&
=\sup _{\eta >0}\bigg|\int_{({\mathbb{R}}^n)^m}K_{u, \eta}(x,y_1,...,y_m)\prod_{j=1}^m(b_j(x)-b_j(y_j))\prod_{i=1}^mf_{i}(y_{i})d\vec{y}\bigg|
\endaligned
$$
and
$$\aligned
V_{\Pi b}^*(\vec{f})(x))&=
\sup_{\eta>0}\bigg|[b_1,[b_2,\cdots[b_{m-1},[b_m,V_\eta ]_m]_{m-1}\cdots]_2]_1 (\vec{f})(x)\bigg|
\\&=\sup _{\eta >0}\bigg|\int_{({\mathbb{R}}^n)^m}K_{v, \eta}(x,y_1,...,y_m)\prod_{j=1}^m(b_j(x)-b_j(y_j))\prod_{i=1}^mf_{i}(y_{i})d\vec{y}\bigg|.
\endaligned
$$

Following \cite{PPTT}, for positive integers $m$ and $j$ with $1\leq
j \leq m$, we denote by $C_j^m$ the family of all finite subsets
$\sigma=\{\sigma(1),\cdots,\sigma(j)\}$ of $\{1,\cdots,m\}$ of $j$
different elements, where we always take $\sigma(k)<\sigma(j)$ if
$k<j$. For any $\sigma \in C_j^m$, we associated the complementary
sequence $\sigma' \in C_j^{m-j}$ given by
$\sigma'=\{1,\cdots,m\}\backslash{\sigma}$ with the convention
$C_0^m=\emptyset$. Given an m-tuple of functions $b$ and $\sigma \in
C_j^m$, we also use the notation $b_{\sigma}$ for the $j$-tuple
obtained from $b$ given by $(b_{\sigma(1)},\cdots,b_{\sigma(j)})$.
\\

Similarly to the above definition for $U_{\Pi b}^*(\vec{f})(x)$ and
$U_{\Pi b}^*(\vec{f})(x)$, $\sigma\in C_j^m$, and
$b_{\sigma}=(b_{\sigma(1)},\cdots,b_{\sigma(j)})$ in $BMO^j$, the
iterated commutator
$$\aligned
U_{\Pi b_{\sigma}}^*(\vec{f})(x)=\sup _{\eta
>0}\bigg|\int_{({\mathbb{R}}^n)^m}K_{u,
\eta}(x,y_1,...,y_m)\prod_{i=1}^j(b_{\sigma(i)}(x)-b_{\sigma(i)}(y_{\sigma(i)}))\prod_{i=1}^mf_{i}(y_{i})d\vec{y}\bigg|;
\endaligned
$$
$$\aligned
V_{\Pi b_{\sigma}}^*(\vec{f})(x)=\sup _{\eta
>0}\bigg|\int_{({\mathbb{R}}^n)^m}K_{v,
\eta}(x,y_1,...,y_m)\prod_{i=1}^j(b_{\sigma(i)}(x)-b_{\sigma(i)}(y_{\sigma(i)}))\prod_{i=1}^mf_{i}(y_{i})d\vec{y}\bigg|;
\endaligned
$$
$$\aligned
I_{\alpha,\Pi
b_{\sigma}}(\vec{f})(x)&=\int_{{(\mathbb{R}^n)}^m}\frac{1}{{|(x-y_1,
\cdots, x-y_m)|}^{mn - \alpha}}\prod
_{i=1}^j(b_{\sigma(i)}(x)-b_{\sigma(i)}(y_{\sigma(i)})) \prod_{i =
1}^{m}f_i(y_i)\,d\vec{y}.
\endaligned
$$

While $\sigma=\{j\}$, $U_{\Pi b_{\sigma}}^*(\vec{f})=
U^*_{b_j}(\vec{f})$,$V_{\Pi
b_{\sigma}}^*(\vec{f})=V^*_{b_j}(\vec{f})$ and $I_{\alpha,\Pi
b_{\sigma}}(\vec{f})=I_{\vec{b},\alpha}^j(\vec{f}).$ If $\sigma=\{1,...,m\}$, then
 $U_{\Pi b_{\sigma}}^*(\vec{f})=
U^*_{\Pi b}(\vec{f})$,$V_{\Pi
b_{\sigma}}^*(\vec{f})=V^*_{\Pi b}(\vec{f})$ and $I_{\alpha,\Pi
b_{\sigma}}(\vec{f})=I_{\alpha,\Pi b}^j(\vec{f}).$

\begin{proposition}[Pointwise control of $M_{\delta}^{\sharp}(U_{\Pi b}^*(\vec{f})), M_{\delta}^{\sharp}(V_{\Pi b}^*(\vec{f})), M_{\delta}^{\sharp}(I_{\alpha, \Pi b}(\vec{f}))$]
Let $0 < \delta < \varepsilon$, $0 < \delta < \frac{1}{m}$ and $0 <
\alpha< mn$. Then there is a constant $C > 0$ depending on $\delta$
and $\varepsilon$ such that
\begin{equation}\label{commutator
pointwise estimate1}
\begin{aligned} M_{\delta}^{\sharp}(U_{\Pi b}^*(\vec{f}))(x)
\leqslant C {\prod_{j=1}^m\|{b_j}\|_{BMO}} \big(\mathcal{M}_{L(\log
L)}(\vec{f})(x) + M_{\varepsilon}(U^*(\vec{f}))(x)\big)
\\ \quad +C\sum_{j=1}^{m-1}\sum_{\sigma\in C_j^m}\prod_{i=1}^j
\|b_{\sigma(i)}\|_{BMO}M_{\varepsilon}(U_{\Pi
b_{\sigma'}}^*(\vec{f}))(x),
\end{aligned}
\end{equation}
\begin{equation}\label{commutator
pointwise estimate2}
\begin{aligned}M_{\delta}^{\sharp}(V_{\Pi b}^*(\vec{f}))(x)
\leqslant C {\prod_{j=1}^m\|{b_j}\|_{BMO}} \big(\mathcal{M}_{L(\log
L)}(\vec{f})(x) + M_{\varepsilon}(V^*(\vec{f}))(x)\big)\\ \quad
+C\sum_{j=1}^{m-1}\sum_{\sigma\in C_j^m}\prod_{i=1}^j
\|b_{\sigma(i)}\|_{BMO}M_{\varepsilon}(V_{\Pi
b_{\sigma'}}^*(\vec{f}))(x),
\end{aligned}
\end{equation}

\begin{equation}\label{commutator
pointwise estimate3}
\begin{aligned}
M_{\delta}^{\sharp}(I_{\alpha, \Pi b}(\vec{f}))(x) \leqslant C
{\prod_{j=1}^m\|{b_j}\|_{BMO}} \big(\mathcal{M}_{L(\log L),
\alpha}(\vec{f})(x)+ M_{\varepsilon}(I_{\alpha}(\vec{f}))(x)\big)\\
\quad +C\sum_{j=1}^{m-1}\sum_{\sigma\in C_j^m}\prod_{i=1}^j
\|b_{\sigma(i)}\|_{BMO}M_{\varepsilon}(I_{\alpha,\Pi
b_{\sigma'}}(\vec{f}))(x).\end{aligned}
\end{equation}
(3.3) still hold for
$\delta=1/m.$
\end{proposition}
{\bf Proof of Proposition 3.1.}

We only give the proof for $U_{\Pi b}^*(\vec{f})$ and $I_{\alpha, \Pi b}(\vec{f})$, since the proof for $V_{\Pi b}^*(\vec{f})$ is almost the same as $U_{\Pi b}^*(\vec{f})$.

For simplicity, we only prove for the case $m=2$, since there is no essential difference for the general case. Fix $b_1,b_2 \in BMO$ and denote any constants by $\rho_1,\, \rho_2$. We split $U_{\Pi b}^*(\vec{f})(x)$ in the following way,
$$\aligned
U_{\Pi b}^*(\vec{f})(x)&=\sup_{\eta>0}|(b_1(x)-\rho_1)(b_2(x)-\rho_2)U_{\eta}
(\vec{f})(x)-(b_1(x)-\rho_1)U_{\eta}(f_1,(b_2-\rho_2)f_2)(x)\\&\quad -(b_2(x)
-\rho_2)U_{\eta}((b_1-\rho_1)f_1,f_2)(x)+U_{\eta}((b_1-\rho_1)f_1,(b_2-\rho_2)f_2)(x)|
\\&\quad =\sup_{\eta>0}|-(b_1(x)-\rho_1)(b_2(x)-\rho_2)U_{\eta}
(\vec{f})(x)+(b_1(x)-\rho_1)U_{\eta,b_2-\rho_2}^2(f_1,f_2)(x)\\&\quad
+(b_2(x)
-\rho_2)U_{\eta,b_1-\rho_1}^1(f_1,f_2)(x)+U_{\eta}((b_1-\rho_1)f_1,(b_2-\rho_2)f_2)(x)|.
\endaligned
$$
Here we denote $U_{\eta,b_1-\rho_1}^1(f_1,f_2)(x)=U_{\eta}((b_1-\rho_1)f_1,f_2)(x)$ and $U_{\eta,b_2-\rho_2}^2(f_1,f_2)(x)=U_{\eta}(f_1,(b_2-\rho_2)f_2)(x)$,
similar notation will be used in the rest of this paper.

Fix $x_0 \in \mathbb{R}^n$ and let $Q$ be a cube centered at $x_0$.
Since $0 < \delta < \frac{1}{m}$, Let $c=\sup_\eta|\sum_{j=1}^3c_j|$, then we have
$${\bigg(\frac{1}{|Q|}\int_Q\big|{|U_{\Pi b}^*(\vec{f})(z)|}^\delta -
{|c|}^\delta\big|\,dz\bigg)}^{\frac{1}{\delta}}  \leqslant  C (T_1 +
T_2+T_3+T_4),$$ where $$T_1 = {\bigg(\frac{1}{|Q|}\int_Q{\big|(b_1(z)-\rho_1)(b_2(z)-\rho_2)\big|^{\delta}U^*(\vec{f})(z)}^{\delta}\,dz\bigg)}^{\frac{1}{\delta}},$$
$$
T_2 = {\bigg(\frac{1}{|Q|}\int_Q{\sup_{\eta>0}\big|(b_1(z)-\rho_1)
[U_{\eta,b_2-\rho_2}^2(f_1,f_2)(z)]\big|}^{\delta}\,dz\bigg)}^{\frac{1}{\delta}}.$$
$$
T_3 = {\bigg(\frac{1}{|Q|}\int_Q{\sup_{\eta>0}\big|(b_2(z)
-\rho_2)[U_{\eta,b_1-\rho_1}^1(f_1,f_2)(z)]\big|}^{\delta}\,dz\bigg)}^{\frac{1}{\delta}}$$
and $$ T_4 =
{\bigg(\frac{1}{|Q|}\int_Q{\sup_{\eta>0}\big|U_{\eta}((b_1-\rho_1)f_1,(b_2-\rho_2)f_2)(z)
-\sum_{j=1}^3c_j\big|}^{\delta}\,dz\bigg)}^{\frac{1}{\delta}}$$ Let
$\rho_j=(b_j)_{3Q}$ be the average of $b_j$ on $3Q$ for $j=1,2.$

For any $1<r_1,r_2,r_3<\infty$ with $\frac 1{r_1}+\frac 1{r_2}+\frac 1{r_3}=1$ and
$r_3 < \frac{\varepsilon}{\delta}$, $T_1$ can be estimated
by using the Holder's inequality and
(\ref{s4}).$$\aligned
T_1& \le C {\bigg(\frac{1}{|Q|}\int_Q{\big|b_1(z) -
\rho_1\big|}^{\delta r_1}\,dz\bigg)}^{\frac{1}{\delta r_1}}{\bigg(\frac{1}{|Q|}\int_Q{\big|b_2(z) -
\rho_2\big|}^{\delta r_2}\,dz\bigg)}^{\frac{1}{\delta r_2}}\\&\quad\times
{\bigg(\frac{1}{|Q|}\int_Q{\big|U^*(\vec{f})(z)\big|}^{\delta
r_3}\,dz\bigg)}^{\frac{1}{\delta r_3}}\\
& \le  C
{\prod_{j=1}^{2}\|b_j\|_{BMO}}M_{\varepsilon}(U^*(\vec{f}))(x_0).\endaligned
$$
Since $T_2$ and $T_3$ are symmetric we only estimate $T_2$. Let
$1<t_1, t_2<\infty$ with $1=1/t_1+1/t_2$ and
$t_2<\frac{\varepsilon}{\delta}$, then $T_1$ can be estimated by
using the H\"{o}lder's inequality and Jensen's inequalities,
$$\aligned
T_2& \le C {\bigg(\frac{1}{|Q|}\int_Q{\big|b_1(z) -
\rho_1\big|}^{\delta t_1}\,dz\bigg)}^{\frac{1}{\delta
t_1}}{\bigg(\frac{1}{|Q|}\int_Q{\sup_{\eta>0}\big|U_{\eta,b_2-\rho_2}^2(f_1,f_2)(z)\big|}^{\delta
t_2}\,dz\bigg)}^{\frac{1}{\delta t_2}}\\
& \le  C\|b_1\|_{BMO}
M_{\varepsilon}(U^{*,2}_{b_2-\rho_2}(\vec{f}))(x_0)
\\
& \le  C \|b_1\|_{BMO} M_{\varepsilon}(U^{*,2}_{b_2}(\vec{f}))(x_0).
\endaligned
$$
Similarly,
$$\aligned
T_3& \le  C\|b_2\|_{BMO}
M_{\varepsilon}(U^*_{b_1-\rho_1}(\vec{f}))(x_0) & \le  C
\|b_2\|_{BMO} M_{\varepsilon}(U^*_{b_1}(\vec{f}))(x_0).
\endaligned
$$

For $T_4$, we denote that $f_i^0 = f_i\chi_{3Q}$ and $f_i^\infty =
f_i - f_i^0$. Note that $c=\sup_\eta|\sum_{j=1}^3c_j|$, where

$$\aligned
c_1=U_{\eta}((b_1-\rho_1)f_1^0,(b_2-\rho_2)f_2^\infty)(x_0),
\endaligned$$

$$\aligned
c_2=U_{\eta}((b_1-\rho_1)f_1^\infty,(b_2-\rho_2)f_2^0)(x_0),
\endaligned$$

$$\aligned
c_3=U_{\eta}((b_1-\rho_1)f_1^\infty,(b_2-\rho_2)f_2^\infty)(x_0).
\endaligned$$

we may split it in the following way
$$\aligned
T_4\le T_{4,1}+T_{4,2}+T_{4,3}+T_{4,4},
\endaligned
$$
where
$$T_{4,1}={\bigg(\frac{1}{|Q|}\int_Q{\sup_{\eta>0}\big|U_{\eta}((b_1-\rho_1)f_1^0,(b_2-\rho_2)f_2^0)(z) \big|}^{\delta}\,dz\bigg)}^{\frac{1}{\delta}},$$
$$T_{4,2}={\bigg(\frac{1}{|Q|}\int_Q{\sup_{\eta>0}\big|U_{\eta}((b_1-\rho_1)f_1^0,(b_2-\rho_2)f_2^\infty)(z) -U_{\eta}((b_1-\rho_1)f_1^0,(b_2-\rho_2)f_2^\infty)(x_0)
\big|}^{\delta}\,dz\bigg)}^{\frac{1}{\delta}},$$
$$T_{4,3}={\bigg(\frac{1}{|Q|}\int_Q{\sup_{\eta>0}\big|U_{\eta}((b_1-\rho_1)f_1^\infty,(b_2-\rho_2)f_2^0)(z) -U_{\eta}((b_1-\rho_1)f_1^\infty,(b_2-\rho_2)f_2^0)(x_0)
\big|}^{\delta}\,dz\bigg)}^{\frac{1}{\delta}}$$ and
$$T_{4,4}={\bigg(\frac{1}{|Q|}\int_Q{\sup_{\eta>0}\big|U_{\eta}((b_1-\rho_1)f_1^\infty,(b_2-\rho_2)f_2^\infty)(z) -U_{\eta}((b_1-\rho_1)f_1^\infty,(b_2-\rho_2)f_2^\infty)(x_0)
\big|}^{\delta}\,dz\bigg)}^{\frac{1}{\delta}}.$$

We consider the first term. Use the Kolmogorov's inequality, lemma
2.2 (a), Theorem B with $w_i\equiv 1$ for $m=2$ and (2.6), then we
deduce that$$\aligned
T_{4,1}&\le C\bigg(\frac{1}{|Q|}\int_{Q}{\Big|U^*((b_1- \rho_1)f_1^0, (b_2- \rho_2)f_2^{0})(z)\Big|}^{{p_0\delta}}\,dz\bigg)^{1/p_0\delta}\\
&\le C {|Q|}^{-2}
{\big\|U^*((b_1- \rho_1)f_1^0, (b_2 - \rho_2)f_2^{0})\big\|}_{L^{\frac{1}{2 }, \infty}(Q)}\\
& \le C {|Q|}^{-2} {\|(b_1 - \rho_1)\|f_1^0\|}_{L^1(Q)}{\|(b_2 - \rho_2)f_2^{0}\|}_{L^1(Q)}\\
& \le C {\|b_1\|}_{BMO}{\|f_1^{0}\|}_{L(\log
L)}{\|b_2\|}_{BMO}{\|f_2^{0}\|}_{L(\log L)}\\
& \le C  \prod_{i=1}^m{\|b_i\|}_{BMO}\mathcal{M}_{L(\log
L)}(\vec{f})(x_0).
\endaligned
$$
By mean value theorem we deduce
$$\aligned
T_{4,2}&\le \frac {C}{|Q|}\int_Q  \sup_{\eta>0}\bigg|U_{\eta}((b_1-
\rho_1)f_1^0, (b_2- \rho_2)f_2^{\infty})(z)- U_{\eta}((b_1-
\rho_1)f_1^{0},(b_2- \rho_2)f_2^{\infty})(x_0)\bigg|dz\\& \le C
\frac 1{|Q|}\int_{Q}\int_{3Q}|(b_1-
\rho_1)f_1(y_1)|dy_1\int_{(3Q)^c}\frac
{|x_0-z|^\varepsilon|b_2(y_2)-\rho_2||f_2(y_2)|dy_2}{(|z-y_1|+|z-y_2|)^{2n+\varepsilon}}dz\\&\le
C\sum _{j=1}^\infty \frac
{j|Q|^{\varepsilon/n}}{(3^j|Q|^{1/n})^{2n+\varepsilon}}\int_{3^{j+1Q}}|(b_1-
\rho_1)f_1(y_1)|dy_1\int_{3^{j+1Q}}|b_2(y_2)-\rho_2||f_2(y_2)|dy_2\\&\le
C\sum_{j=1}^\infty
\frac{1}{3^{j\varepsilon}}\prod_{i=1}^2{\|b_i\|}_{BMO}\|f_i\|_{L(\log
L),3^{j+1Q}}
\\&\le
C \prod_{i=1}^2{\|b_i\|}_{BMO}\mathcal{M}_{L(\log L)}(\vec{f})(x_0).
\endaligned
$$
Similarly as $T_{4,2}$, we can get the estimates for $T_{4,3}$. Now
we are in a position to deal $T_{4,4}$. Note
that\begin{eqnarray*}\lefteqn{\big||U_{\eta}((b_1-
\rho_1)f_1^{\infty}, (b_2- \rho_2)f_2^{\infty})(z)
 - {(|U_{\eta}(((b_1-
\rho_1)f_1^{\infty},(b_2- \rho_2)f_2^{\infty}))(x_0)}\big|}\\
& \leqslant &C
\int_{{(\mathbb{R}^n\backslash3Q)}^2}\frac{{|Q|}^{\frac{\varepsilon}{n}}|(b_1-
\rho_1)||b_2(y_2) - \rho_2|}{{|(x_0 - y_1, x_0 - y_2)|}^{2n+
\varepsilon}}\prod_{i =
1}^2|f_i^\infty(z_i)|\,d\vec{y}\\
 & \leqslant &C
\sum_{k = 1}^\infty\int_{{(3^{k + 1}Q)}^2 \backslash
{(3^kQ)}^2}\frac{{|Q|}^{\frac{\varepsilon}{n}}|(b_1-
\rho_1)||b_2(y_2) - \rho_2|}{{(3^k{|Q|}^{\frac{1}{n}})}^{2n +
\varepsilon}}\prod_{i =
1}^2|f_i^\infty(y_i)|\,d\vec{y}\\
& \leqslant & C \prod_{i=1}^2{\|b_i\|}_{BMO}\mathcal{M}_{L(\log
L)}(\vec{f})(x_0).\end{eqnarray*} Thus, we have $$T_{4,4}\le C
\prod_{i=1}^2{\|b_i\|}_{BMO}\mathcal{M}_{L(\log L)}(\vec{f})(x_0).$$
Thus we complete the proof of this lemma for $U_{\Pi b}^*(\vec
{f})$.


Next, we prove (3.3) for $I_{{\alpha,\Pi b}}(\vec{f})$, we split
$$\aligned
I_{\alpha, \Pi
b}(\vec{f})(x)&=(b_1(x)-\rho_1)(b_2(x)-\rho_2)I_{\alpha}(\vec{f})(x)-(b_1(x)-\rho_1)
I_{\alpha}(f_1,(b_2-\rho_2)f_2)(x)\\&\quad -(b_2(x)
-\rho_2)I_{\alpha}((b_1-\rho_1)f_1,f_2)(x)+I_{\alpha}((b_1-\rho_1)f_1,(b_2-\rho_2)f_2)(x)
\\&=-(b_1(x)-\rho_1)(b_2(x)-\rho_2)I_{\alpha}(\vec{f})(x)+(b_1(x)-\rho_1)
I_{b_2-\rho_2,\alpha}^2(f_1,f_2)(x)\\&\quad +(b_2(x)
-\rho_2)I_{b_1-\rho_1,\alpha}^1(f_1,f_2)(x)+I_{\alpha}((b_1-\rho_1)f_1,(b_2-\rho_2)f_2)(x).
\endaligned
$$
Fix $x_0 \in \mathbb{R}^n$ and let $Q$ be a cube centered at $x_0$.
Denote any constants by $c=(I_{\alpha}(f_1^{0},(b_2-\rho_2)f_2^{\infty})(x_0) +
I_{\alpha}(f_1^{\infty},
(b_2-\rho_2)f_2^{0})(x_0)+I_{\alpha}(f_1^{\infty},
(b_2-\rho_2)f_2^{\infty})(x_0))=:c_1+c_2+c_3.$ Then

Since $0 < \delta \leqslant \frac{1}{m}$, then we have
$${\bigg(\frac{1}{|Q|}\int_Q\big|{|I_{\alpha, \Pi b}(\vec{f})(z)|}^\delta -
{|c|}^\delta\big|\,dz\bigg)}^{\frac{1}{\delta}}  \leqslant  C (S_1 +
S_2+S_3+S_4),$$ where
$$S_1 = {\bigg(\frac{1}{|Q|}\int_Q{\big|(b_1(z)-\rho_1)(b_2(z)-\rho_2)\big|^\delta
\big|I_{\alpha}(\vec{f})(z)\big|}^{\delta}\,dz\bigg)}^{\frac{1}{\delta}},$$
$$
S_2 = {\bigg(\frac{1}{|Q|}\int_Q{\big|(b_1(x)-\rho_1)
I_{b_2-\rho_2,\alpha}^2(f_1,f_2)(z)
\big|}^{\delta}\,dz\bigg)}^{\frac{1}{\delta}}.$$
$$
S_3 = {\bigg(\frac{1}{|Q|}\int_Q{\big|(b_2(x)
-\rho_2)I_{b_1-\rho_1,\alpha}^1(f_1,f_2)(z)
\big|}^{\delta}\,dz\bigg)}^{\frac{1}{\delta}}$$ and $$ S_4 =
{\bigg(\frac{1}{|Q|}\int_Q{\big|I_{\alpha}((b_1-\rho_1)f_1,(b_2-\rho_2)f_2)(z)
- c\big|}^{\delta}\,dz\bigg)}^{\frac{1}{\delta}}.$$ Let
$\rho_j=(b_j)_{3Q}$ be the average of $b_j$ on $3Q$ for $j=1,2.$

For any $1<r_1,r_2,r_3<\infty$ with $\frac 1{r_1}+\frac 1{r_2}+\frac 1{r_3}=1$ and $r_3 < \frac{\varepsilon}{\delta}$, $S_1$ can be estimated
by using the Holder's inequality and
(\ref{s4}).$$\aligned
S_1& \le C {\bigg(\frac{1}{|Q|}\int_Q{\big|b_1(z) -
\rho_1\big|}^{\delta r_1}\,dz\bigg)}^{\frac{1}{\delta r_1}}{\bigg(\frac{1}{|Q|}\int_Q{\big|b_2(z) -
\rho_2\big|}^{\delta r_2}\,dz\bigg)}^{\frac{1}{\delta r_2}}\\&\quad\times
{\bigg(\frac{1}{|Q|}\int_Q{\big|I_{\alpha}(\vec{f})(z)\big|}^{\delta
r_3}\,dz\bigg)}^{\frac{1}{\delta r_3}}\\
& \le  C
{\prod_{j=1}^{2}\|b_j\|_{BMO}}M_{\varepsilon}(I_\alpha(\vec{f}))(x_0).\endaligned
$$

As the argument of $T_2$, we still take $1<t_1, t_2<\infty$ with
$1=1/t_1+1/t_2$ and $t_2<\frac{\varepsilon}{\delta}$
$$\aligned
S_2 &= {\bigg(\frac{1}{|Q|}\int_Q{\big|(b_1(x)-\rho_1)
I_{b_2-\rho_2,\alpha}^2(f_1,f_2)(z)
\big|}^{\delta}\,dz\bigg)}^{\frac{1}{\delta}}\\
& \le  C
\|b_1\|_{BMO}\mathcal{M}_{t_2\delta}(I_{b_2-\rho_2,\alpha}^2(f_1,f_2))(x_0)\\
& \le  C
\|b_1\|_{BMO}\mathcal{M}_{\varepsilon}(I_{b_2-\rho_2,\alpha}^2(f_1,f_2))(x_0).
\endaligned
$$

Similarly, we can get the estimates for $S_3$ as we deal $S_2$.
Next, for $S_4$, we denote that $f_i^0 = f_i\chi_{3Q}$ and
$f_i^\infty = f_i - f_i^0$ and Let
$c=(I_{\alpha}(f_1^{0},(b_2-\rho_2)f_2^{\infty})(x_0) +
I_{\alpha}(f_1^{\infty},
(b_2-\rho_2)f_2^{0})(x_0)+I_{\alpha}(f_1^{\infty},
(b_2-\rho_2)f_2^{\infty})(x_0))$, then $S_4$ can be written as
$$\aligned
S_4\le S_{4,1}+S_{4,2}+S_{4,3}+S_{4,4},
\endaligned
$$
where
$$S_{4,1}={\bigg(\frac{1}{|Q|}\int_Q{\big|I_{\alpha}((b_1-\rho_1)f_1^0,(b_2-\rho_2)f_2^0)(z) \big|}^{\delta}\,dz\bigg)}^{\frac{1}{\delta}},$$
$$S_{4,2}={\bigg(\frac{1}{|Q|}\int_Q{\big|I_{\alpha}((b_1-\rho_1)f_1^0,(b_2-\rho_2)f_2^\infty)(z) -I_{\alpha}((b_1-\rho_1)f_1^0,(b_2-\rho_2)f_2^\infty)(x_0)
\big|}^{\delta}\,dz\bigg)}^{\frac{1}{\delta}},$$
$$S_{4,3}={\bigg(\frac{1}{|Q|}\int_Q{\big|I_{\alpha}((b_1-\rho_1)f_1^\infty,(b_2-\rho_2)f_2^0)(z) -I_{\alpha}((b_1-\rho_1)f_1^\infty,(b_2-\rho_2)f_2^0)(x_0)
\big|}^{\delta}\,dz\bigg)}^{\frac{1}{\delta}}$$
and
$$S_{4,4}={\bigg(\frac{1}{|Q|}\int_Q{\big|I_{\alpha}((b_1-\rho_1)f_1^\infty,(b_2-\rho_2)f_2^\infty)(z) -I_{\alpha}((b_1-\rho_1)f_1^\infty,(b_2-\rho_2)f_2^\infty)(x_0)
\big|}^{\delta}\,dz\bigg)}^{\frac{1}{\delta}}.$$

Use H\"{o}lder inequality, the Kolmogorov's inequality
(\ref{kolmogorov}) when $p = \frac{1}{2}$ and $q = \frac{n}{2n -
\alpha}$, (2.11) in Lemma 2.4, then we deduce that$$\aligned
S_{4,1}&\le C\bigg(\frac{1}{|Q|}\int_{Q}{\Big|I_\alpha((b_1- \rho_1)f_1^0, (b_2- \rho_2)f_2^{0})(z)\Big|}^{\frac{1}{2}}\,dz\bigg)^2\\
&\le C {|Q|}^{\frac{\alpha}{n}-2}
{\big\|I_\alpha((b_1- \rho_1)f_1^0, (b_2 - \rho_2)f_2^{0})\big\|}_{L^{\frac{n}{2n - \alpha}, \infty}(Q)}\\
& \le C {|Q|}^{\frac{\alpha}{n}-2} {\|(b_1 - \rho_1)f_1^{0}\|}_{L^1(Q)}{\|\|(b_2 - \rho_2)f_2^{0}\|\|}_{L^1(Q)}\\
& \le C
{|3Q|}^{\frac{\alpha}{n}}{\|b_1\|}_{BMO}{\|f_1^{0}\|}_{L(\log L),
Q}{\|b_2\|}_{BMO}{\|f_2^{0}\|}_{L(\log L),
Q}\\
& \le C  \prod_{j=1}^2{\|b_j\|}_{BMO}\mathcal{M}_{L(\log L),
\alpha}(\vec{f})(x_0).
\endaligned
$$
By mean value theorem again, we deduce
$$\aligned
S_{4,2}&\le \frac {C}{|Q|}\int_Q  \bigg|I_\alpha((b_1- \rho_1)f_1^0,
(b_2- \rho_2)f_2^{\infty})(z) -I_{\alpha}((b_1- \rho_1)f_1^{0},(b_2-
\rho_2)f_2^{\infty})(x_0)\bigg|dz\\& \le C \frac
1{|Q|}\int_{3Q}|(b_1(y_1)- \rho_1)f_1(y_1)|dy_1\int_{(3Q)^c}\frac
{|x_0-z||(b_1-
\rho_1)||b_2(y_2)-\rho_2||f_2(y_2)|dy_2}{(|z_1-y_1|+|z_2-y_2|)^{2n-\alpha+1}}dz\\&\le
C\sum _{j=1}^\infty \frac
j{(3^j|Q|^{1/n})^{2n-\alpha+1}}\int_{3Q}|(b_1(y_1)-
\rho_1)f_1(y_1)|dy_1\int_{3^{j+1Q}}|b_2(y_2)-\rho_2||f_2(y_2)|dy_2\\&\le
C \prod_{j=1}^2{\|b_j\|}_{BMO}\mathcal{M}_{L(\log L),
\alpha}(\vec{f})(x_0).
\endaligned
$$
Similarly as $S_{4,2}$, we can get the estimates for $S_{4,3}$. Now
we are in a position to deal $S_{4,4}$.
\begin{eqnarray*}\lefteqn{\big|I_\alpha((b_1- \rho_1)f_1^{\infty}, (b_2- \rho_2)f_2^{\infty})(z)
 - {(I_\alpha((b_1- \rho_1)f_1^{\infty},(b_2- \rho_2)f_2^{\infty}))(x_0)}\big|}\\
& \leqslant &C
\int_{{(\mathbb{R}^n\backslash3Q)}^2}\frac{{|Q|}^{\frac{1}{n}}|b_1(y_1)
- \rho_1||b_2(y_2) - \rho_2|}{{|(x_0 - y_1, x_0 - y_2)|}^{2n -
\alpha + 1}}\prod_{i =
1}^2|f_i^\infty(y_i)|\,d\vec{y}\\
 & \leqslant &C
\sum_{k = 1}^\infty\int_{{(3^{k + 1}Q)}^2 \backslash
{(3^kQ)}^2}\frac{{|Q|}^{\frac{1}{n}}|b_1(y_1) - \rho_1||b_2(y_2) -
\rho_2|}{{(3^k{|Q|}^{\frac{1}{n}})}^{2n - \alpha + 1}}\prod_{i =
1}^2|f_i^\infty(y_i)|\,d\vec{y}\\
& \leqslant & C \sum_{k =
1}^\infty\frac{k}{3^k}{\|b_2\|}_{BMO}{|3^{k +
1}Q|}^{\frac{\alpha}{n}}\prod_{j=1}^2{\|f_j^\infty\|}_{L(\log L),
3^{k +
1}Q}\\
& \leqslant & C \prod_{j=1}^2{\|b_j\|}_{BMO}\mathcal{M}_{L(\log L),
\alpha} (\vec{f})(x_0).\end{eqnarray*} So we obtain
 $$S_{4,i}\le C {\|b_1\|}_{BMO}{\|b_2\|}_{BMO}\mathcal{M}_{L(\log L), \alpha}
(\vec{f})(x_0).$$ Thus we complete the proof for this lemma.
\begin{proposition}[Pointwise control of $M_{\delta}^{\sharp}(U^*(\vec{f})), M_{\delta}^{\sharp}(V^*(\vec{f})), M_{\delta}^{\sharp}(I_{ \alpha}(\vec{f}))$]
Let $0 < \delta < \varepsilon$, $0 < \delta <\frac{1}{m}$
and $0 < \alpha <m n$. Then
there is $C > 0$ depending on $\delta$ and $\varepsilon$ such that
\begin{equation}\label{commutator
pointwise estimate}M_{\delta}^{\sharp}(U^*(\vec{f}))(x) \leqslant C\mathcal{M}(f)(x),\end{equation}
\begin{equation}\label{commutator
pointwise estimate}M_{\delta}^{\sharp}(V^*(\vec{f}))(x) \leqslant C \mathcal{M}(f)(x),\end{equation}
\begin{equation}\label{commutator
pointwise estimate}M_{1/m}^{\sharp}(I_{\alpha}(\vec{f}))(x) \leqslant  \mathcal{M}_\alpha(f)(x).\end{equation} for all bounded $\vec{f}$ with compact support.
\end{proposition}

{\bf Proof.}

The proof of (3.4) and (3.5) follows from similar steps in Theorem 3.2 of \cite{new maximal function} and combine the method we used in the above proposition, here we omit the proof. On the other hand, (2.7) has already been obtained in \cite{CX}, Proposition 5.2.

Now, we can obtain
\begin{theorem} Let $0<p$ and $w\in A_\infty$. Suppose that $\vec{b} \in (BMO)^m$. Then there is a constant $C$ independent of $\vec{b}$ and a constant $C_1 $ (may dependent on $\vec{b}$) such that
\begin{equation}\int_{\Bbb {R}^n}|U_{\Pi b}^*(\vec{f})(x)|^p\omega(x)dx \le C\prod_{i=1}^m\|b_i\|_{BMO}\int_{\Bbb {R}^n}[\mathcal{M}_{L(\log L)}(f)(x)]^pw(x)dx,\end{equation}
\begin{equation}\aligned \sup_{t>0} \frac{1}{\Phi^m(1/t)}w(\{y\in \Bbb {R}^n: |U_{\Pi b}^*&\vec{f}(y)|>t^m\})\\ &\le C_1\sup \frac{1}{\Phi^m(1/t)}w(\{y\in \Bbb {R}^n: \mathcal{M}_{L(\log L)}(f)(y)>t^m\}).\endaligned \end{equation}
 Similar results hold for $V_{\Pi b}^*(\vec{f})$.

\end{theorem}
{\bf Proof of Theorem 3.1.}The proof of the above Theorem 3.1 are
now standard as the case for multilinear C-Z singular integral
operators. We briefly indicate such arguments in the case m=2, but,
as the reader will immediately notice, and iterative procedure using
(3.1) and (3.2)can be followed to obtain the general case. \\
Using Fefferman-Stein inequality and pointwise estimate in
proposition 3.1 we will have

$$\aligned
\|U_{\Pi b}^*(\vec{f})\|_{L^p(\omega)}&\leq\|M_{\delta}(U_{\Pi
b}^*(\vec{f}))\|_{L^p(\omega)}\leq C \|M^{\sharp}_{\delta}(U_{\Pi b}^*(\vec{f}))\|_{L^p(\omega)}\\
&\leq C \prod_{i=1}^2
\|b_i\|_{BMO}\biggl(\|\mathcal{M}_{L(logL)}(\vec{f})\|_{L^p(\omega)}+\|M^{\sharp}_{\varepsilon}(U^*(\vec{f}))\|_{L^p(\omega)}
\biggl) \\ & \quad
+C\biggl(\|b_2\|_{BMO}\|M^{\sharp}_{\varepsilon}(U_{b_1}^{*}(\vec{f}))\|_{L^p(\omega)}+
\|b_1\|_{BMO}\|M^{\sharp}_{\varepsilon}(U_{b_2}^{*}(\vec{f}))\|_{L^p(\omega)}\biggl).
\endaligned
$$
Hence, next we estimate
$\|M^{\sharp}_{\varepsilon}(U_{b_2}^{*}(\vec{f}))\|_{L^p(\omega)}$,
$\|M^{\sharp}_{\varepsilon}(U_{b_1}^{*}(\vec{f}))\|_{L^p(\omega)}$ has
the similar estimate.
Set $c_{\eta}=U_{\eta}(f_1^{0},(b_2- \rho_2)f_2^{\infty})(x_0) +
U_{\eta}(f_1^{\infty}, (b_2-
\rho_2)f_2^{0})(x_0)+U_{\eta}(f_1^{\infty}, (b_2-
\rho_2)f_2^{\infty})(x_0)$ and $c=\sup_{\eta>0}\{|c_{\eta}|\}$, then
$$\aligned
|U_{b_2}^*(\vec{f})(z)-c|&\leq \sup _{\eta
>0}\bigg|\int_{({\mathbb{R}}^n)^2}K_{u,
\eta}(z,y_1,y_2)((b_2(z)-\rho_2)-(b_2(y_2)-\rho_2))\prod_{i=1}^2f_{i}(y_{i})d\vec{y}+c_{\eta}\bigg|\\
&\leq C |b_2(z)-\rho_2|U^*(f_1,f_2)(z)+\sup _{\eta
>0}|U_{\eta}(f_1,(b_2-\rho_2)f_2)(z)-c_{\eta}|.
\endaligned
$$

For arbitrary $0<\varepsilon'<\frac{1}{2}$, take $1<t_1, t_2<\infty$
with $1=1/t_1+1/t_2$ and $t_2<\frac{\varepsilon'}{\varepsilon}$, we
have
$$\aligned
&{}{\bigg(\frac{1}{|Q|}\int_Q \big|(b_2(z)-\rho_2)U^*(f_1,f_2)(z)
\biggl|^{\varepsilon}\,dz\bigg)}^{\frac{1}{\varepsilon}}\\&\le
{\bigg(\frac{1}{|Q|}\int_Q |b_2(z)-\rho_2|
^{t_1\varepsilon}\,dz\bigg)}^{\frac{1}{t_1\varepsilon}}
\bigg(\frac{1}{|Q|}\big|\int_Q U^*(f_1,f_2)(z)
\big|^{t_2\varepsilon}\,dz\bigg)^{\frac{1}{t_2\varepsilon}}
\\&\le C\|b_2\|_{BMO}\mathcal{M}_{\varepsilon'}(U^*(f_1,f_2))(x_0)
.\endaligned
$$

As the proof of Proposition 3.1, then $U_{\eta}(f_1, (b_2-
\rho_2)f_2)$ can be written as
$$\aligned
U_{\eta}(f_1, (b_2- \rho_2)f_2) & =U_{\eta}(f_1^{0}, (b_2
-\rho_2)f_2^{0})
 + U_{\eta}(f_1^{0},(b_2- \rho_2)f_2^{\infty})\\
& \quad +  U_{\eta}(f_1^{\infty}, (b_2-
\rho_2)f_2^{0})+U_{\eta}(f_1^{\infty}, (b_2 -
\rho_2)f_2^{\infty}).\endaligned
$$
Take $1<p_0<1/(2\varepsilon)$ and using H\"{o}lder's inequality
again, we have
$$\aligned
&{}{\bigg(\frac{1}{|Q|}\int_Q{\sup_{\eta>0}\big|U_{\eta}(f_1,(b_2-\rho_2)f_2)(z)
-
c_{\eta}\big|}^{\varepsilon}\,dz\bigg)}^{\frac{1}{\varepsilon}}\\&\le
\bigg(\frac
1{|Q|}\int_Q\sup_{\eta>0}\bigg|U_{\eta}(f_1,(b_2-\rho_2)f_2)(z)-c_{\eta}\bigg|^{p_0\varepsilon}
dz\bigg)^{1/p_0\varepsilon}\\&\le ( G_{1}+G_{2}+G_{3}+G_{4}),
\endaligned
$$
where
$$G_{1}=\bigg(\frac{1}{|Q|}\int_{Q}{\sup_{\eta>0}\Big|U_{\eta}(f_1^0,
(b_2 -
\rho_2)f_2^{0})(z)\Big|}^{p_0\varepsilon}\,dz\bigg)^{1/{p_0}\varepsilon},$$
$$G_{2}=\bigg(\frac{1}{|Q|}\int_{Q}{\sup_{\eta>0}\Big|U_{\eta}(f_1^0, (b_2 - \rho_2)f_2^{\infty})(z)
-U_{\eta}(f_1^{0},(b_2- \rho_2)f_2^{\infty})(x_0)
\Big|}^{p_0\varepsilon}\,dz\bigg)^{1/{p_0}\varepsilon},$$
$$G_{3}=\bigg(\frac{1}{|Q|}\int_{Q}{\sup_{\eta>0}\Big|U_{\eta}(f_1^\infty,
(b_2 - \rho_2)f_2^{0})(z)-U_{\eta}(f_1^{\infty}, (b_2 -
\rho_2)f_2^{0})(x_0)\Big|}^{p_0\varepsilon}\,dz\bigg)^{1/{p_0}\varepsilon}$$
and
$$G_{4}=\bigg(\frac{1}{|Q|}\int_{Q}{\sup_{\eta>0}\Big|U_{\eta}(f_1^\infty,
 (b_2- \rho_2)f_2^{\infty})(z)-U_{\eta}(f_1^{\infty},
(b_2-
\rho_2)f_2^{\infty})(x_0)\Big|}^{p_0\varepsilon}\,dz\bigg)^{1/{p_0}\varepsilon}.$$
The similar procedure for $T_4$ in the Proposition 3.1, we obtain
$$\aligned
G_{1}\le C  {\|b_2\|}_{BMO}\mathcal{M}_{L(\log L)}^2(\vec{f})(x_0).
\endaligned
$$
By mean value theorem we deduce
$$\aligned
G_{2}\le C {\|b_2\|}_{BMO}\mathcal{M}_{L(\log L)}^2(\vec{f})(x_0).
\endaligned
$$
Similarly as $G_{2}$, we can get the estimates for $G_{3}$. Moreover
$$G_{4}\le C {\|b_2\|}_{BMO}\mathcal{M}_{L(\log L)}^2
(\vec{f})(x_0).$$ By proposition 3.2, so we have
$$\aligned
\|M^{\sharp}_{\varepsilon}[U_{b_2}^*(\vec{f})]\|_{L^p(\omega)}& \leq
C\|b_2\|_{BMO}(\|\mathcal{M}(\vec{f})\|_{L^p(\omega)}+\|\mathcal{M}_{L(\log
L)}^2(\vec{f})\|_{L^p(\omega)})\\& \leq
C\|b_2\|_{BMO}\|\mathcal{M}_{L(\log L)}(\vec{f})\|_{L^p(\omega)}.
\endaligned
$$

The desired inequality now follows. Since the left main steps and
the ideas are almost the same as \cite{PPTT}, here we omit the
proof. So we get the estimate of strong type and weak type.

{\bf Proof of Theorem 1.1-1.2.} Theorem 1.1 follows by the reason
that $T_{*,\Pi b}(\vec{f})\le U_{\Pi b}^*(\vec{f})(x)+V_{\Pi
b}^*(\vec{f})(x)$, Theorem 3.1 and the weighted strong boundedness
of $\mathcal{M}_{L(\log L)}$ in \cite{new maximal function}. Theorem
1.2 follows by repeating the same steps as in \cite{new maximal
function}, \cite{PPTT} and the method used in \cite{zhang}. Since
the main steps and the ideas are almost the same, here we omit the
proof.

{\bf Proof of Theorem 1.3.} Theorem 1.3 follows by using Proposition
3.1 and the estimate for $I_{\vec{b},\alpha}^j (j=1,2)$, which is
Theorem 2.7 in \cite{CX}.

\section[Weighted end-point estimates for $I_{\alpha, \Pi b}(\vec{f})$]{Weighted end-point estimates for $I_{\alpha, \Pi b}(\vec{f})$}\label{commutators 2}

Firstly, we will consider the end-point estimate of multilinear fractional
$L(\log L)$ type maximal operator.

\begin{proposition}[Weighted end-point estimate for $\mathcal{M}_{L(\log L), \alpha}$]
Let $\Phi(t) = t(1 + \log^+t)$ and $\vec{\omega} \in A_{((1,
\cdots, 1), \frac{n}{mn - \alpha})}$.
If $0<\alpha<mn$, then there is a $C > 0$ such
that
\begin{equation}\aligned
\label{end unweight comm}{\nu_{\vec{\omega}}}&^{\frac{n}{mn - \alpha}}\bigg(\Big\{x\in \mathbb{R}^n
: \mathcal{M}_{L(\log L),\alpha}(\vec{f})(x) > t^{\frac{mn -
\alpha}{n}}\Big\}\bigg) \\& \leqslant C
 {\bigg\{ \bigg[1 + \frac{\alpha}{mn}\log^+\bigg(\prod_{i = 1}^m\int_{\mathbb{R}^n}\Phi^{(m)}\Big(\frac{|f_i(y_i)|}{t}\Big)\,dy_i\bigg)\bigg]^m \prod_{j = 1}^m\int_{\mathbb{R}^n}\Phi^{(m)}\Big(\frac{|f_j(y_j)|}{t}\Big)\omega_j(y_j)\,dy_j\bigg\}}^{\frac{n}{mn - \alpha}}.\endaligned
\end{equation}
If $0<\alpha_j<n$ for each $1\le j\le m$, $\sum_{j=1}^m\alpha_j=\alpha,$ then there is a $C > 0$ such
that
\begin{equation}\aligned
\label{end unweight comm}{\nu_{\vec{\omega}}}&^{\frac{n}{mn - \alpha}}\bigg(\Big\{x\in \mathbb{R}^n
: \mathcal{M}_{L(\log L),\alpha}(\vec{f})(x) > t^{\frac{mn -
\alpha}{n}}\Big\}\bigg) \\& \leqslant C
 {\bigg\{\prod_{j = 1}^m\bigg[1 + \frac{\alpha_j}{n}\log^+\bigg(\prod_{i = 1}^m\int_{\mathbb{R}^n}\Phi^{(m)}\Big(\frac{|f_i(y_i)|}{t}\Big)\,dy_i\bigg)\bigg]\int_{\mathbb{R}^n}\Phi^{(m)}
 \Big(\frac{|f_j(y_j)|}{t}\Big)\omega_j(y_j)\,dy_j\bigg\}}^{\frac{n}{mn - \alpha}}.\endaligned
\end{equation}

\end{proposition}

\textbf{Proof.}
By the homogeneity, we can assume $t = 1$. We first prove (4.2). Denote
that
\begin{center}$E_{1} = \Big\{x\in \mathbb{R}^n : \mathcal{M}_{L(\log L), \alpha}(\vec{f})(x) >
1\Big\}$ and $E_{1, k} = E_{1} \cap B(0, k)$,\end{center}where $B(0,
k) = \{x\in \mathbb{R}^n : |x| \leqslant k\}$. By the monotone
convergence theorem, it suffices to estimate $E_{1, k}$.

For any $x \in E_{1, k}$, there is a cube $Q_x$ such that
\begin{equation}\label{bigger that 1}1 < {|Q_{x}|}^{\frac{\alpha}{n}}\prod_{j=1}^m{\|f_j\|}_{L(\log L), Q}
\, .\end{equation}Hence,
${\{Q_x\}}_{x \in E_{1, k}}$ is a family of cubes covering $E_{1,
k}$. Using a covering argument, we obtain a finite family of disjoint
cubes $\{Q_{x_l}\}$ whose dilations cover $F$ such that
\begin{equation}\displaystyle |E_{1, k} |\le C\sum_l|
Q_{x_l}| \quad \hbox{and} \quad 1 < {|Q_{x_l}|}^{\frac{\alpha}{n}}\prod_{j=1}^m{\|f_j\|}_{L(\log L), Q_{x_l}}
\, .\end{equation}

We follow the main steps first as in \cite {PPTT} and denote $C_h^m$ to be the family of all subset $\sigma=(\sigma(1),...,\sigma(h))$ from $\{1,...,m\}$ with $1\le h\le m$ different elements. Given $\sigma\in C_h^m$ and a cube $Q_{x_l}$, if $|Q_{x_l}|^{\alpha_{\sigma(j)}}{\|f_{\sigma(j)}\|_{L(\log L), Q_{x_l}}}{}>1$ for $j=1,...,h$, we say that $j\in B_\sigma$ and $|Q_{x_l}|^{\alpha_{\sigma(j)}}{\|f_{\sigma(j)}\|_{L(\log L), Q_{x_l}}}{}\le 1$ for $j=h+1,...,m$. Denote $$A_k=\prod_{j=1}^k{|Q_{x_l}|^{\alpha_{\sigma(j)}/n}\|f_{\sigma(j)}\|}_{L(\log L), Q_{x_l}}$$ and $A_0=1$. Then it is easy to check that if $\sigma\in C_h^m$ and  $j\in B_\sigma$, for any $1\le k\le m$, we have $A_k>1$ and
$$\aligned
1&<\prod_{j=1}^k{|Q_{x_l}|^{\alpha_{\sigma(j)}/n}\|f_{\sigma(j)}\|}_{L(\log L), Q_{x_l}} \\& ={\bigg\||Q_{x_l}|^{\alpha_{\sigma(k)}/n}f_{\sigma(k)}A_{k-1}\bigg\|_{\Phi,Q_{x_l}}}.\endaligned$$ Or, equivalently
\begin{equation}1<\frac {1}{{|Q_{x_l}|} } \int_{Q_{x_l}}\Phi\bigg(|Q_{x_l}|^{\alpha_{\sigma(k)}/n}f_{\sigma(k)}\bigg(\prod_{j=1}^{k-1}{|Q_{x_l}|^{\alpha_{\sigma(j)}/n}\|f_{\sigma(j)}\|}_{L(\log L), Q_{x_l}}\bigg)\bigg).\end{equation}
By the following equivalence
$$\|f\|_{\Phi,Q}\simeq \inf_{\mu>0}\{\mu+\frac {\mu}{|Q_{x_l}|}\int_{Q_{x_l}} \Phi(|f|/\mu)\}.$$
If $1\le j\le m-h-1,$ we obtain
$$\aligned
\Phi^j(A_{m-j})&=\Phi^j(\||Q_{x_l}|^{{\alpha_{\sigma(m-j)}/n}}f_{\sigma(m-j)}A_{m-j-1}\|_{\Phi, Q_{x_l}}).
\endaligned
$$
Since$\||Q_{x_l}|^{{\alpha_{\sigma(m-j)}/n}}f_{\sigma(m-j)}A_{m-j-1}\|_{\Phi, Q}>1$, Using the fact that $\Phi$
is submultiplicative (i.e. $\Phi(st) \leqslant \Phi(s)\Phi(t)$ for
$s, t > 0$) and Jensen's inequality, we have
$$\aligned
\Phi^j(A_{m-j})&=\Phi^j(\||Q_{x_l}|^{{\alpha_{\sigma(m-j)}/n}}f_{\sigma(m-j)}A_{m-j-1}\|_{\Phi, Q})\\&
\le C\Phi^j(1+\frac 1 {|Q_{x_l}|}\int_Q\Phi(|Q_{x_l}|^{{\alpha_{\sigma(m-j)}/n}}f_{\sigma(m-j)}A_{m-j-1}))\\& \le C\frac 1 {|Q_{x_l}|}\int_Q\Phi^{j+1}(|Q_{x_l}|^{{\alpha_{\sigma(m-j)}/n}}f_{\sigma(m-j)})\Phi^{j+1}(A_{m-j-1}).
\endaligned
$$
By iterating the inequalities above and the fact that $\||Q_{x_l}|^{\alpha_{\sigma(j)}/n}
f_{\sigma(j)}\|_{\Phi,Q_{x_l}}>1$ for $j\in B_\sigma$, $\Phi^{j+1}\le \Phi^m$ and $\Phi^{m-h+1}\le \Phi^m$ for $1\le h\le m$ and $0\le j\le m-h-1$, we have
\begin{equation}\aligned
1&<\frac {1}{{|Q_{x_l}|} } \int_{Q_{x_l}}\Phi\big(|Q_{x_l}|^{\alpha_{\sigma(m)}/n}f_{\sigma(m)}\big)\frac {1}{|Q_{x_l}|} \int_{Q_{x_l}}\Phi^2(|Q_{x_l}|^{{\alpha_{\sigma(m-1)}/n}}f_{\sigma(m-1)})\Phi^2({A_{m-2}})\\& \le
\bigg( \prod_{j=0}^{m-h-1}\frac {1}{|Q_{x_l}|}\int_{Q_{x_l}}\Phi^{j+1}(|Q_{x_l}|^{{\alpha_{\sigma(m-1)}/n}}f_{\sigma(m-1)})\bigg)\prod_{j=1}^h\Phi^{m-h}(\||Q_{x_l}|^{\alpha_{\sigma(j)}/n}
f_{\sigma(j)}\|_{\Phi,Q_{x_l}})\\& \le \bigg( \prod_{j=0}^{m-h-1}\frac {1}{|Q_{x_l}|}\int_{Q_{x_l}}\Phi^{j+1}(|Q_{x_l}|^{{\alpha_{\sigma(m-1)}/n}}f_{\sigma(m-1)})\bigg)\prod_{j=1}^h  \frac {1}{|Q_{x_l}|}\int_{Q_{x_l}} \Phi^{m-h+1}(|Q_{x_l}|^{\alpha_{\sigma(j)}/n}
f_{\sigma(j)})\\& \le C\prod_{j=1}^m
\frac {1}{|Q_{x_l}|}\int_{Q_{x_l}}\Phi^{m}(|Q_{x_l}|^{\frac {\alpha_j}{n}}f_j).
\endaligned
\end{equation}
We obtain
\begin{equation}\aligned
1&< C\prod_{j=1}^m
\frac {1}{|Q_{x_l}|}\int_{Q_{x_l}}\Phi^{m}(|Q_{x_l}|^{\frac {\alpha_j}{n}}f_j)\\&\le C
\prod_{j=1}^m
\frac {1}{|Q_{x_l}|}\int_{Q_{x_l}}\Phi^{m}(|Q_{x_l}|^{\frac {\alpha_j}{n}})\Phi^m(f_j)\\&\le C
\prod_{j=1}^m\frac {1}{|Q_{x_l}|} |Q_{x_l}|^{\frac {\alpha_j}{n}}\big(1+\log^+ |Q_{x_l}|^{{\frac {\alpha_j}{n}}}\big)\int_{Q_{x_l}}\Phi^{m}(f_j).
\endaligned
\end{equation}
Since $\alpha_j<n, $ there exists a constant $C_0>1$ and $\eta_1, ..., \eta_m$ small enough, such that
$$0<\eta_j<1-\frac {\alpha_j}n,\quad \quad 1+\log^+ t^{{\frac {\alpha_j}{n}}}\le t^{\eta_j}\quad \hbox{\ \ if}\quad t>C_0.$$
Denote $\eta=\sum_{j=1}^m \eta_j$, then by (4.7) if $|Q_{x_l}|>C_0$ we have
\begin{equation}
|Q_{x_l}|^{m-\frac \alpha n-\eta}\le C\prod_{j=1}^m
\int_{Q_{x_l}}\Phi^{m}(f_j).
\end{equation}
Thus,
$$(m-\frac{ \alpha}{ n}-\eta) \log^+(|Q_{x_l}|^{\frac {\alpha_j}{n}})\le C \frac {\alpha_j}{n}\log^+ \bigg(\prod_{j=1}^m
\int_{Q_{x_l}}\Phi^{m}(f_j)\bigg).
$$
By (4.7) again, we have
\begin{equation}\aligned
|Q_{x_l}|^{m-\frac{ \alpha}{ n}}&\le C \prod_{j=1}^m \bigg\{1+\frac {\alpha_j}{n}\log^+ \bigg(\prod_{j=1}^m
\int_{Q_{x_l}}\Phi^{m}(f_j)\bigg)\bigg\}
\int_{Q_{x_l}}\Phi^{m}(f_j).
\endaligned
\end{equation}
On the other hand, if $|Q_{x_l}|\le C_0$, then it is easy to see $1+\log^+ |Q_{x_l}|^{{\frac {\alpha_j}{n}}}\le C.$ Thus
\begin{equation}|Q_{x_l}|^{m-\frac{ \alpha}{ n}}\le C \prod_{j=1}^m
\int_{Q_{x_l}}\Phi^{m}(f_j).\end{equation}
(4.9) and (4.10) yield that
\begin{equation}\aligned
|Q_{x_l}|^{m-\frac{ \alpha}{ n}}&\le C \prod_{j=1}^m \bigg\{1+\frac {\alpha_j}{n}\log^+ \bigg(\prod_{j=1}^m
\int_{Q_{x_l}}\Phi^{m}(f_j)\bigg)\bigg\}
\int_{Q_{x_l}}\Phi^{m}(f_j).
\endaligned
\end{equation}
Finally, by (4.4) and the definition of class $A_{((1, \cdots,
1), \frac{n}{mn - \alpha})}$, we
$$ \aligned
{\bigg(\int_{E_{1, k}}{\nu_{\vec{\omega}}}^{\frac{n}{m n - \alpha}}\,\bigg)}^{\frac{m n -
\alpha}{n}}& \le { \bigg(
\sum_{h=1}^m \sum_{\sigma\in C_h^m} \sum_{l\in B_\sigma} \int_{Q_{x_l}} {\nu_{\vec{\omega}}}^{\frac{n}{m n -
\alpha}}\,\bigg)}^{\frac{m n -
\alpha}{n}}\\&\le C{
\sum_{h=1}^m \sum_{\sigma\in C_h^m} \sum_{l\in B_\sigma}\bigg(\int_{Q_{x_l}} {\nu_{\vec{\omega}}}^{\frac{n}{m n -
\alpha}}\,\bigg)}^{\frac{m n -
\alpha}{n}}\\&\le
C\sum_{h=1}^m \sum_{\sigma\in C_h^m} \sum_{l\in B_\sigma}|Q_{x_l}|^{m-\frac {\alpha}{n}}\prod_{j=1}^m\inf w_j\\&\le
C \sum_{h=1}^m \sum_{\sigma\in C_h^m} \sum_{l\in B_\sigma} \prod_{j=1}^m \bigg\{1+\frac {\alpha_j}{n}\log^+ \bigg(\prod_{j=1}^m
\int_{Q_{x_l}}\Phi^{m}(f_j)\bigg)\bigg\}
\int_{Q_{x_l}}\Phi^{m}(f_j) w_j\\&\le C\sum_{h=1}^m \sum_{\sigma\in C_h^m} \sum_{l\in B_\sigma} \prod_{j=1}^m \bigg\{1+\frac {\alpha_j}{n}\log^+ \bigg(\prod_{j=1}^m
\int_{\mathbb{R}^n}\Phi^{m}(f_j)\bigg)\bigg\}
\int_{Q_{x_l}}\Phi^{m}(f_j) w_j\\&\le C\prod_{j=1}^m \bigg\{1+\frac {\alpha_j}{n}\log^+ \bigg(\prod_{j=1}^m
\int_{\mathbb{R}^n}\Phi^{m}(f_j)\bigg)\bigg\}
\int_{\mathbb{R}^n}\Phi^{m}(f_j) w_j.
\endaligned
$$
The proof of inequality (4.2) is finished.

Inequality (4.1) follows by taking $\alpha_j=\alpha/m<n$ in the above proof.

\textbf{Proof of Theorem 1.4 and Corollary 1.1.}

To prove Theorem 1.4, we follow the main steps as in \cite{CX}, without changes till the last step by using the above Proposition 3.1, We will obtain Theorem 1.4.

To prove Corollary 1.1, similarly as in linear case \cite{DLZ}, we define \begin{equation}\aligned
\overline{I_{\alpha, \Pi b}}(\vec{f})(x)&=\int_{{(\mathbb{R}^n)}^m}\frac{\prod _{j=1}^m|b_j(x)-b_j(y_j)|}{{|(x-y_1, \cdots, x-y_m)|}^{mn - \alpha}}
\prod_{i = 1}^{m}|f_i(y_i)|\,d\vec{y},
\endaligned
\end{equation}
And careful check in the proof of Theorem 1.3-1.4 shows that Theorem 1.3-1.4 still hold for $\overline{I_{\alpha, \Pi b}}.$ Note the fact that $\mathcal{M}_{{\Pi b}, \alpha}(\vec{f})(x)\le \overline{I_{\alpha, \Pi b}}(|f_1|,...,|f_m|)(x)$, this implies Corollary 1.1.
\subsection*{Acknowledgements} The author want to express his sincerely thanks to the unknown referee for his or her valuable remarks
and pointing out a problem in the previous version of this manuscript.

\begin{flushleft}

\vspace{0.3cm}\textsc{Qingying Xue\\School of Mathematical
Sciences\\Beijing Normal University\\Laboratory of Mathematics and
Complex Systems\\Ministry of Education\\Beijing 100875\\
People's
Republic of China\\
}

\emph{E-mail address}: \textsf{qyxue@bnu.edu.cn}
\end{flushleft}


\begin{thebibliography}{20}

\bibitem{Alva}J. Alvarez, \emph{Continuity properties for linear commutators of
Calder\'{o}n-Zygmund operators}, Collect. Math. \textbf{49} (1998),
no. {1}, 17-3.

\bibitem{BGP}A. Bernardis, O. Gorosito, G. Pradolini, \emph {Weighted inequalities for multilinear potential operators and its commutators}, preprint, arXiv:1007.0445.
   \bibitem{C}X. Chen,
\emph{Weighted estimates for maximal Operator of
multilinear singular integral}, to appear in Bull. Polish Acad. Sci. Math.
\bibitem{CX}X. Chen and Q. Xue,
\emph{Weighted estimates for a class of multilinear fractional type
operators}, J. Math. Anal. Appl. \textbf{362} (2010), no. 2,
 355-373.
\bibitem{Polynomial}M. Christ and J.-L. Journ\'{e},
\emph{Polynomial growth estimates for multilinear singular integral operators},
 Acta Math. \textbf{159} (1987), 51-80.

 \bibitem{on commutators}R. R. Coifman and Y. Meyer,
\emph{On commutators of singular integrals and bilinear singular integrals},
Trans. Amer. Math. Soc. \textbf{212} (1975), 315-331.

\bibitem{operators}R. R. Coifman and Y. Meyer,
\emph{Commutateurs d'integrales singulires et op\'{e}rateurs
multilin\'{e}aires}, Ann. Inst. Fourier (Grenoble) \textbf{28} (1978),
No. 3, 177-202.
\bibitem{CRW}R. R. Coifman, R. Rochberg, Guido Weiss, \emph {Factorization theorems for Hardy spaces in several variables} Ann. of Math. (2) \textbf{103} (1976), 611-635.


\bibitem{DTT}C. Demeter, T. Tao and C. Thiele, \emph{Maximal multilinear operators}, Trans. Amer. Math. Soc.
\textbf{360} (2008), 4989-5042.

\bibitem{DLZ} Y. Ding, S. Lu and P. Zhang, \emph {Weak estimates for commutators of fractional integral
operators}, Science in China (Series A), English Series
2001, \textbf{44} (2002), No.7, 877-888.

\bibitem{weighted norm inequ}J. Garc\'{i}a-Cuerva, J. L.
Rubio de Francia, \emph{Weighted Norm Inequalities and Related
Topics}, North-Holland Math. Studies \textbf{116}, North-Holland,
Amsterdam, 1985.
\bibitem{multi frac grafakos}L. Grafakos, \emph{On multilinear fractional
integrals}, Studia Math. \textbf{102}(1992), 49-56.

\bibitem{GLPT}L. Grafakos, L. Liu, C. P\'{e}rez and R.H. Torres, \emph{The multilinear strong maximal function}, preprint.
\bibitem{on Hardy spaces}L. Grafakos and N. Kalton,
\emph{ Multilinear Calder\'{o}n-Zygmund operators on Hardy spaces},
Collect. Math. \textbf{52} (2001), No. {2}, 169-179.
\bibitem{multi C-Z}L. Grafakos and R. H. Torres, \emph{Multilinear Calder\'{o}n-Zygmund
theory}, Adv. Math. \textbf{165} (2002), No. {1}, 124-164.
\bibitem{report}L. Grafakos and R. H. Torres, \emph{On multilinear
singular integrals of Calder\'{o}n-Zygmund type}, Proceedings of the
6th International Conference on Harmonic Analysis and Partial
Differential Equations (El Escorial). Publ. Mat. 2002, Vol. Extra,
57-91.
\bibitem{GT1}L. Grafakos and R. H. Torres,
\emph{Maximal operator and weighted norm inequalities for
multilinear singular integrals}, Indiana. Univ. Math. J. \textbf{51}
(2002), No. {5}, 1261-1276.

\bibitem{kenig frac}C. E. Kenig and E. M. Stein, \emph{Multilinear estimates
and fractional integration}, Math. Res. Lett. \textbf{6} (1999),
1-15.
\bibitem{new maximal function}A. K. Lerner, S. Ombrosi, C. P\'{e}rez, R. H.
Torres and R. Trujillo-Gonz\'{a}lez, \emph{New maximal functions and
multiple weights for the multilinear Calder\'{o}n-Zygmund theory},
Adv. Math. \textbf{220} (2009), no. {4}, 1222--1264.

\bibitem{LXY}W. Li, Q. Xue and
K. Yabuta, \emph{Multilinear Calder\'{o}n-Zygmund operators on
weighted Hardy spaces}, Studia Math. \textbf{199} (2010), no. {1}, 1-16.
\bibitem{LXY1}W. Li, Q. Xue and
K. Yabuta, \emph{Maximal operator for multilinear Calder\'{o}n-Zygmund singular integral operators on weighted Hardy spaces}, to appear in J. Math. Anal. Appl. \textbf{373} (2011), no. 2, 384-392.

\bibitem{moen}K. Moen, \emph{Weighted inequalities for multilinear
fractional integral operators}, Collect. Math. \textbf{60}(2009),
213-238.

\bibitem{MW}B. Muckenhoupt and R. Wheeden, \emph{Weighted norm inequalities for  fractional
integral}, Trans. Amer. Math. Soc.  \textbf{192} (1974), 261-274.

\bibitem{pe}C. P\'{e}rez, \emph{end-point estimates for commutators of singular integral operators},
J. Funct. Anal. \text{128} (1995), no. 1, 163¨C185.
\textbf{192}(1974), 261-274.
\bibitem{PPTT}C. P\'{e}rez, G. Pradolini, R.H. Torres, and R. Trujillo-Gonz¨¢lez, \emph{End-point estimates for iterated commutators of multilinear singular integrals}, preprint, arXiv:1004.4976.

\bibitem{pe1}C. P\'{e}rez, R. Trujillo-Gonz¨¢lez,  \emph{Sharp weighted estimates for multilinear commutators}, J. London Math. Soc. (2) \text{65}
(2002), 672¨C692.
\bibitem{orlicz}M. M. Rao and Z. D. Ren,
\emph{Theory of Orlicz Spaces}, Monographs and Textbooks in Pure and
Applied Math. \textbf{146}, Marcel Dekker, New York, 1991.

\bibitem{weight frac welland}G. V. Welland, \emph{Weighted norm
inequalities for fractional integrals}, Proc. Amer. Math. Soc.
\textbf{51}(1975), no. {1}, 143-148.
\bibitem{zhang}P. Zhang, \emph{Weighted estimates for maximal multilinear commutators}, Math. Nachr. \textbf {279}, no. 4, (2006)445-462.

\end{thebibliography}
\end{document}